\documentclass[11pt, a4paper]{article}
\usepackage[T1]{fontenc}
\usepackage{amsfonts}
\usepackage{amsmath}
\usepackage{amssymb}
\usepackage{amsthm}
\usepackage{bbm}
\usepackage{bm}
\usepackage{mathrsfs}
\usepackage{verbatim}
\usepackage{setspace}
\usepackage{color}
\usepackage{pdfsync}
\usepackage{enumitem}

\usepackage[pdfborder={0 0 0}]{hyperref}
\hypersetup{
  urlcolor = black,
  pdfauthor = {Bruno Bouchard, Marcel Nutz},
  pdfkeywords = {Stochastic target; Stochastic differential game; Knightian uncertainty; Shaking of coefficients; Viscosity solution},
  pdftitle = {Stochastic Target Games and Dynamic Programming via Regularized Viscosity Solutions},
  pdfsubject = {Stochastic Target Games and Dynamic Programming via Regularized Viscosity Solutions},
  pdfpagemode = UseNone
}

\theoremstyle{plain}
\newtheorem{theorem}{Theorem}[section]
\newtheorem{proposition}[theorem]{Proposition}
\newtheorem{lemma}[theorem]{Lemma}
\newtheorem{corollary}[theorem]{Corollary}

\theoremstyle{definition}

\newtheorem{remark}[theorem]{Remark}

\newtheorem{assumption}{Assumption}

\theoremstyle{remark}

\renewenvironment{thebibliography}[1]{%
\begin{oldthebibliography}{#1}%
\setlength{\baselineskip}{.98em}
\linespread{1}
\small
\setlength{\parskip}{0.28ex}%
\setlength{\itemsep}{.38em}%
}%
{%
\end{oldthebibliography}%
}

\newcommand{\eps}{\varepsilon}
\newcommand{\vp}{\varphi}

\newcommand{\M}{\mathbb{M}}

\newcommand{\Q}{\mathbb{Q}}
\newcommand{\R}{\mathbb{R}}

\newcommand{\E}{\mathbb{E}}
\newcommand{\F}{\mathbb{F}}

\renewcommand{\P}{\mathbb{P}}

\newcommand{\mfU}{\mathfrak{U}}
\newcommand{\mfu}{\mathfrak{u}}

\def\Ac{{\cal A}}
\def\Bc{{\cal B}}

\def\Dc{{\cal D}}
\def\Ec{{\cal E}}
\def\Fc{{\cal F}}

\def\Lc{{\cal L}}

\def\Oc{{\cal O}}
\def\Pc{{\cal P}}

\def\Tc{{\cal T}}
\def\Uc{{\cal U}}

\def\Yc{{\cal Y}}

\def \E{\mathbb{E}}
\def \F{\mathbb{F}}

\def \M{\mathbb{M}}

\def \R{\mathbb{R}}

\def\P{\mathbb{P}}
\def\Q{\mathbb{Q}}

\def\Sb{{\mathbf S}}
\def\Hb{{\mathbf H}}

\def\Db{{\mathbf D}}

\newcommand{\as}{\mbox{-a.s.}}
\newcommand{\1}{\mathbf{1}}

\def \Int{\displaystyle\int}

\def\esssup{{\rm ess}\!\sup\limits}

\def\ep{\hbox{ }\hfill$\Box$}
\def\reff#1{{\rm(\ref{#1})}}
\def\be{\begin{eqnarray}}
\def\ee{\end{eqnarray}}
\def\beq{\begin{equation}}
\def\eeq{\end{equation}}

\def\proof{{\noindent \bf Proof. }}

\def\as{{\rm a.s.}}

\def\Esp#1{\mathbb{E}\left[#1\right]}

\def\x{\times}

\def\eqref#1{\reff{#1}}
\def\wu{{\underline w}}

\def\D{{\mathcal D}}
\def\Di{{\D_{<T}}}
\def\DT{{\D_{T}}}
\def \[{[\,\!\![}
\def \]{]\,\!\!]}

\numberwithin{equation}{section}

\begin{document}

\title{\vspace{-0em}
Stochastic Target Games and Dynamic Programming via Regularized Viscosity Solutions%
\thanks{We are grateful to Erhan Bayraktar, Pierre Cardaliaguet and two anonymous referees for valuable discussions and suggestions.}
\date{This version January 2015}
\author{
     Bruno Bouchard%
     \thanks{CEREMADE, Universit\'e Paris Dauphine and CREST-ENSAE, bouchard@ceremade.dauphine.fr. Research supported by ANR Liquirisk and Labex ECODEC.
     }
     \and
  Marcel Nutz%
  \thanks{Departments of Statistics and Mathematics, Columbia University, New York, mnutz@columbia.edu. Research supported by NSF Grant DMS-1208985.
  }
 }
}
\maketitle \vspace{-0em}

\begin{abstract}
We study a class of stochastic target games where one player tries to find a strategy such that the state process almost-surely reaches a given target, no matter which action is chosen by the opponent. Our main result is a geometric dynamic programming principle which allows us to characterize the value function as the viscosity solution of a non-linear partial differential equation. Because abstract measurable selection arguments cannot be used in this context, the main obstacle is the construction of measurable almost-optimal strategies. We propose a novel approach where smooth supersolutions are used to define almost-optimal strategies of Markovian type, similarly as in verification arguments for classical solutions of Hamilton--Jacobi--Bellman equations. The smooth supersolutions are constructed by an extension of Krylov's method of shaken coefficients. We apply our results to a problem of option pricing under model uncertainty with different interest rates for borrowing and lending.
\end{abstract}

\vspace{3em}

{\small


\newpage

\section{Introduction}\label{sec: introduction}

We study a stochastic (semi) game where we try to find a strategy $\mfu$ such that the controlled state process almost-surely reaches a given target at the time horizon $T$, no matter which control $\alpha$ is chosen by the adverse player. Here $\alpha$ is any predictable process with values in a given bounded set $A\subset \R^{d}$, whereas $\mfu$ is a non-anticipating strategy which associates to each $\alpha$ a predictable process $\mfu[\alpha]$ with values in a given set $U\subset \R^{d}$. More precisely, if $\Ac$ and $\Uc$ denote the collections of predictable processes with values in $A$ and $U$, respectively, then $\mfu$ is a map $\Ac\to \Uc$ such that $\mfu[\alpha]|_{[0,s]}=\mfu[\alpha']|_{[0,s]}$ on the set $\{\alpha|_{[0,s]}=\alpha'|_{[0,s]}\}$, for all $\alpha,\alpha'\in\Ac$ and $s\le T$. We denote by $\mfU$ the collection of all such strategies.

Given an initial condition $(x,y)\in\R^d\times \R$ at time $t$ and $(\mfu,\alpha)\in \mfU\x \Ac$, the $(d+1)$-dimensional state process $(X^{\alpha}_{t,x},Y^{\mfu,\alpha}_{t,x,y})(s)$, $t\leq s\leq T$ is defined as the solution  of the stochastic differential equation
\begin{eqnarray*}
\left\{\begin{array}{rcl}
 dX(s)&=&\mu_{X}(s,X(s),\alpha_{s}) ds +\sigma_{X}(s,X(s),\alpha_{s})dW_{s},\\[.3em]
 dY(s)&=&\mu_{Y}(s,X(s),Y(s),\mfu[\alpha]_{s},\alpha_{s}) ds +\sigma_{Y}(s,X(s),Y(s),\mfu[\alpha]_{s},\alpha_{s})dW_{s},
\end{array}
\right.
\end{eqnarray*}
where $W$ is a Brownian motion and the coefficients satisfy suitable continuity and growth conditions.
Given a measurable function $g$, the value function of the stochastic target game is then given by
\begin{equation}\label{eq: def v intro}
v(t,x)=\inf\big\{y\in \R: \exists \; \mfu \in \mfU\;\mbox{ s.t.}\; Y^{\mfu,\alpha}_{t,x,y}(T)\ge g(X^{\alpha}_{t,x}(T))\;\as\;\forall\;\alpha \in \Ac\big\}.
\end{equation}
That is, $v(t,x)$ is the smallest $y$ from which we can drive the state process into the epigraph of $g$ by using a strategy $\mfu$ to react to the adverse player.
The aim of this paper is to provide a dynamic programming principle for the stochastic target game and to characterize $v$ in terms of a Hamilton--Jacobi--Bellman equation.

In the case where $A$ is a singleton and $\mfu$ is a control, the above is a standard stochastic target problem in the terminology of~\cite{SonerTouzi.02b,SonerTouzi.02a}. In \cite{SonerTouzi.02b}, it is shown that the value function of this target problem satisfies a geometric dynamic programming principle (GDP) and, consequently, is a discontinuous viscosity solution of an associated Hamilton--Jacobi--Bellman equation. The GDP consists of two parts, called GDP1 and GDP2. Roughly, given a family $\{\theta^{\mfu},\mfu\in \mfU\}$ of stopping times with values in $[t,T]$, GDP1 states that
$$
\mbox{if $y>v(t,x)$, then there exists $\mfu\in \mfU$ such hat $Y^{\mfu,a}_{t,y}(\theta^{\mfu})\ge v(\theta^{\mfu},X_{t,x}^{a}(\theta^{\mfu}))$,}
$$
and conversely, GDP2 states that
$$
\mbox{if there exists $\mfu\in \mfU$ such that $Y^{\mfu,a}_{t,y}(\theta^{\mfu})> v(\theta^{\mfu},X_{t,x}^{a}(\theta^{\mfu}))$, then $y\ge v(t,x).$}
$$
The line of argument for GDP1 can be reproduced in the context of games.
However, all previous proofs of GDP2 crucially rely on the construction of almost optimal controls through measurable selections.
It is well known that in the context of games, when $\mfu$ is a strategy as defined above, the possibility of using such a selection is a completely open problem. The main difficulties come from the lack of separability in the space $\mfU$ of strategies and the irregular dependence on the adverse control $\alpha$.

For zero-sum differential games in standard form, see e.g.\ \cite{FlemingSouganidis.89} and~\cite{BuLi08}, there are by now several workarounds for this problem; they rely on approximations and exploit the continuity properties of the cost functions. We also refer to \cite{PhamZhang.12} and \cite{Sirbu.13} for recent developments in different directions, in a setting where both players use controls. In \cite{BouchardMoreauNutz.12}, stochastic target games were analyzed when the target is of controlled-loss type; that is, of the form $\Esp{\ell(X^{\mfu,\alpha}_{t,x,y}(T), Y^{\mfu,\alpha}_{t,x,y}(T))}\ge m$, where $m$ is given and the loss function $\ell$ has certain continuity properties. Again, the regularity of $\ell$ was exploited to circumvent the problem of measurable selection. By contrast, the almost-sure target of the game under consideration is highly discontinuous, which has prevented us from arguing similarly as in the mentioned works.

In this paper, we follow a completely different and novel idea. As a starting point, recall that in the context of standard control problems with a smooth value function, we can sometimes use a \emph{verification argument.}  Here, an optimal control of Markovian (i.e.\ feedback) type is defined explicitly in terms of the derivatives of the value function. It plays the role of the almost-optimal controls mentioned above and renders measurable selection arguments unnecessary. Of course, this procedure requires a smooth value function, which cannot be expected in our context. However, it will turn out that a smooth supersolution with specific properties can be used in a similar spirit. The outline of our argument runs as follows.
\begin{itemize}
 \item[(a)] Show that the value function $v$ satisfies a version of GDP1 above (Theorem~\ref{thm: GDP1}).
 \item[(b)] Deduce from GDP1 that $v$ is a viscosity supersolution of the associated Hamilton--Jacobi--Bellman equation (Theorem~\ref{thm: super sol property}).
 \item[(c)] Regularize $v$ to find a smooth supersolution $w$ which is close to $v$ is a specific sense (Lemma~\ref{lem : smoothing supersol semi lin conv}).
 \item[(d)] Using $w$, construct a strategy of Markovian type that matches the criterion in~\reff{eq: def v intro} when starting slightly above $v(t,x)$ and use this strategy to prove a version of GDP2 (Theorem~\ref{thm: GDP2}).
  \item[(e)] Deduce from GDP2 that $v$ is also a viscosity subsolution of the Hamilton--Jacobi--Bellman equation (Theorem~\ref{thm: sub sol property}).
\end{itemize}
The arguments for (a) and (b) are along the lines of \cite{SonerTouzi.02b}, while at least part of the proof of~(e) follows~\cite{BouchardMoreauNutz.12}. The construction of the smooth supersolution in (c) is based on Krylov's method of \emph{shaking the coefficients} from \cite[Theorem 2.1]{Krylov00} (see also~\cite{BarlesJakobsen.02}), which we extend here to semi-linear equations by considering value functions of controlled forward-backward stochastic differential equations (Theorem~\ref{thm : smoothing supersol semi lin conv gene}). We mention that this method imposes a concavity condition on one of the operators. On the other hand, to the best of our knowledge, our result is the first proof of dynamic programming for stochastic target games with almost-sure target.\footnote{Note added in proof: In the follow-up work \cite{bayraktar2014stochastic}, a similar target game problem is treated by a stochastic Perron's method which does not require a dynamic programming principle.}

Our results can be compared to the second order backward stochastic differential equations of \cite{SonerTouziZhang.2010bsde}, where the authors use  quasi-sure analysis in the  ``weak formulation''. Their setting is more general in the sense that path-dependence is allowed and concavity is not needed. On the other hand, we allow nonlinear dynamics for $X$, while their setting corresponds to the case where $\sigma_{X}(\cdot,a)=a$ and $\mu_{X}=0$.

We apply our results to a problem from mathematical finance where an option is to be super-hedged under model uncertainty. The model uncertainty is addressed in a worst-case fashion, so that it can be modeled by an adverse player as above. More precisely, the drift and volatility of the underlying of the option as well as the two interest rates for borrowing and lending depend on the adverse control~$\alpha$. Various incarnations of the super-hedging problem have been considered in the recent literature; see~\cite{AvellanedaLevyParas.95, BouchardNutz.13, DenisMartini.06, NeufeldNutz.12,Peng.10,   SonerTouziZhang.2010rep, SonerTouziZhang.2010dual} and the references therein. The now standard approach is to use the weak formulation where the uncertainty is modeled by a non-dominated set $\Pc$ of possible laws on path space. The super-hedging property is then required to hold almost-surely under each element of $\Pc$ and the study involves the difficulty of dealing with a non-dominated set of probabilities (quasi-sure analysis). We adopt here a different point of view, where uncertainty is modeled as a game under a single probability measure; namely, the adverse player (``nature'') chooses  drift and volatility while we can react by using a suitable non-anticipating strategy and the superhedging-property is required to hold almost-surely for any control of the adverse player. Our formulation is thus in the spirit of \cite{TalayZheng.02} and \cite{TevzadzeToronjadzeUzunashvili.13} where problems of portfolio management are phrased as stochastic differential games of standard form in the framework of \cite{FlemingSouganidis.89}.

The remainder of this paper is organized as follows. In Section~\ref{sec: stochastic target games}, we formulate the stochastic target game in detail and provide the geometric dynamic programming principle together with the Hamilton--Jacobi--Bellman equation for the value function $v$. The key part of the dynamic programming principle, GDP2, is based on the regularization result, which is developed in Section~\ref{sec: smoothing} in a slightly more general setting. In Section~\ref{sec: examples}, we exemplify our results by the application to super-hedging under model uncertainty.

\section{Geometric dynamic programming for stochastic target games}\label{sec: stochastic target games}

In this section, we first detail our problem formulation. In the second subsection, we provide the first part of the geometric dynamic programming principle, GDP1, and infer the supersolution property of the value function. In the third subsection, we prove the difficult part, GDP2, together with the subsolution property.

\subsection{Problem formulation}\label{sec: problem formulation}

Fix a time horizon $T>0$, let $\Omega$ be the space of continuous functions $\omega: [0,T]\to \R^{d}$ and let $\P$ be the Wiener measure on $\Omega$.
Moreover, let $W$ be the canonical process on $\Omega$, defined by $W_t(\omega)=\omega_t$. We denote by $\F=(\Fc_s)_{0\leq s\leq T}$ the augmented filtration generated by $W$. Furthermore, for $0\leq t\leq T$, we denote by $\F^t=(\Fc^t_s)_{0\leq s\leq T}$ the augmented filtration generated by $(W_s-W_t)_{s\geq t}$; by convention, $\Fc^t_s$ is trivial for $s\leq t$. We denote by $\Uc^{t}$ (resp.\ $\Ac^{t}$) the collection of all $\F^{t}$-predictable processes in $L^p(\P\otimes dt)$ with values in a given Borel subset $U$ (resp.\ bounded {closed} subset $A$) of $\R^{d}$, where $p\geq2$ is fixed throughout.
Finally, let
$$
\D:=[0,T]\x \R^{d},\quad\Di:=[0,T)\x \R^{d},\quad\DT:=\{T\}\x \R^{d}.
$$

Given $(t,x,y)\in \D\x \R$ and $(\nu,\alpha)\in \Uc^{t}\x \Ac^{t}$, we let $(X^{\alpha}_{t,x},Y^{\nu,\alpha}_{t,x,y})(s)$, $t\leq s\leq T$ be the unique strong solution of
\be\label{eq: dynamics}
\left\{\begin{array}{rcl}
dX(s)&=&\mu_{X}(s,X(s),\alpha_{s}) ds +\sigma_{X}(s,X(s),\alpha_{s})dW_{s},\\[.3em]
dY(s)&=&\mu_{Y}(s,X(s),Y(s),\nu_{s}, \alpha_{s}) ds +\sigma_{Y}(s,X(s),Y(s),\nu_{s},\alpha_{s})dW_{s}
\end{array}
\right.
\ee
with initial data $(X(t),Y(t))=(x,y)$. The coefficients $\mu_{X}, \mu_{Y}, \sigma_{Y}$ and $\sigma_{X}$ are supposed to be continuous in all variables, {uniformly in the last one},  and take values in $\R^{d}$, $\R$, $\R^{d}$ and \mbox{$\M^{d}:=\R^{d\times d}$,} respectively. (Elements of $\R^{d}$ as viewed as column vectors).  We assume throughout that there exists $K>0$ such that
\be\label{eq: cond mu sigma}
\begin{array}{c}
|\mu_{X}(\cdot,x,\cdot)-\mu_{X}(\cdot,x',\cdot)|+|\sigma_{X}(\cdot,x,\cdot)-\sigma_{X}(\cdot,x',\cdot)|\le K|x-x'| ,\\[.3em]
|\mu_{X}(\cdot,x,\cdot)|+|\sigma_{X}(\cdot,x,\cdot)|\le  K,\\[.3em]
|\mu_{Y}(\cdot,y, \cdot)-\mu_{Y}(\cdot,y', \cdot)| + |\sigma_{Y}(\cdot,y, \cdot)-\sigma_{Y}(\cdot,y', \cdot)| \le K |y-y'| ,\\[.3em]
 |\mu_{Y}(\cdot,y, u,\cdot)|+|\sigma_{Y}(\cdot,y,u, \cdot)|\le K(1+|u|+|y|)
\end{array}
\ee
for all $(x,y),(x',y')\in \R^d\x\R$ and $u\in U$. In particular, this ensures that the SDE~\reff{eq: dynamics} is well-posed. Moreover, we can note that the solution is in fact adapted not only to $\F$ but also to $\F^t$.

\begin{remark} For simplicity, we consider the above dynamics for any initial point $x\in \R^{d}$. The case where $X^{\alpha}_{t,x}$ is restricted to an open domain $\Oc\subset\R^d$ will be discussed in Remark~\ref{rem: open domain} below. 
%
\end{remark}
For the derivation of the viscosity supersolution property, we shall also impose a condition on the growth of $\mu_Y$ relative to $\sigma_Y$:
\be\label{eq: cond mu over sigma}
 \sup_{u\in U}  \frac{|\mu_{Y}(\cdot, u,\cdot)|}{1+|\sigma_{Y}(\cdot, u,\cdot)|}\;\mbox{ is locally bounded.}
\ee

Let $t\le T$. We say that a map $\mfu:\Ac^{t}\to \Uc^{t}$, $\alpha\mapsto \mfu[\alpha]$ is a $t$-admissible strategy if it is non-anticipating in the sense that
\be\label{eq: cond non anticipativity}
 \{\omega \in \Omega: \alpha(\omega){|_{[t,s]}}= \alpha'(\omega){|_{[t,s]}}\}\subset  \{\omega \in \Omega: \mfu[\alpha](\omega) {|_{[t,s]}}= \mfu[\alpha'](\omega) {|_{[t,s]}}\}\;\as
 \ee
 for all $s\in [t,T]$ and $\alpha,\alpha'\in \Ac^{t}$, where $|_{[t,s]}$ indicates the restriction to the interval $[t,s]$.  We denote by $\mfU^{t}$  the collection of all $t$-admissible strategies; moreover, we write $Y^{\mfu,\alpha}_{t,x,y}$ for $Y^{\mfu[\alpha],\alpha}_{t,x,y}$. Finally, let $g: \R^d\to\R$ be a measurable function; then we can introduce the value function of our stochastic target game,
\begin{equation}\label{eq: def v}
v(t,x):=\inf\big\{y\in \R: \exists \; \mfu \in \mfU^{t}\;\mbox{ s.t.}\; Y^{\mfu,\alpha}_{t,x,y}(T)\ge g(X^{\alpha}_{t,x}(T))\; \as\;\forall\;\alpha \in \Ac^{t}\big\}
\end{equation}
for $(t,x)\in \Dc$. We shall assume throughout that
\be\label{eq: ass v bounded g bounded lipschitz}
\mbox{$g$ is bounded  and Lipschitz continuous  {on $\R^{d}$}, and  $v$ is bounded {on $\cal D$}.}
\ee

\begin{remark} {The condition that $v$ is bounded has to be checked on a case-by-case basis. One typical example in which $v^{+}$ is bounded is when there exists $u$ such that 
$\sigma_{Y}(\cdot,u,\cdot)=0$. Then, the condition \reff{eq: cond mu sigma} on $\mu_{Y}$ implies that $v^{+}$ is bounded by $e^{KT} (KT+\sup_{\R^{d}} g)$. 
A simple situation in which $v^{-}$ is bounded  is when $\sigma_{Y}^{-1}\mu_{Y}$ is bounded. Then, an obvious change of measure argument allows to turn $Y^{\mfu,\alpha}_{t,x,y}$ into a martingale, for $(\mfu,\alpha)$ given, which implies that $v\ge \inf_{\R^{d}} g$. See also Section \ref{sec: examples}. 
}
\end{remark}

\subsection{First part of the dynamic programming principle and supersolution property}

We first provide one side of the geometric dynamic programming principle, GDP1, for the value function $v$ of~\reff{eq: def v}.
We denote by $v_{*}$ the lower-semicontinuous envelope of $v$ on $\D$.

\begin{theorem}[{\bf GDP1}]\label{thm: GDP1} Let $(t,x,y)\in \D\x \R$ and let $\{\theta^{\mfu,\alpha},(\mfu,\alpha)\in \mfU^{t}\x \Ac^{t}\}$ be a family of $\F^{t}$-stopping times with values in $[t,T]$. Assume that $y>v(t,x)$. Then, there exists $\mfu \in \mfU^{t}$ such that
$$
Y_{t,x,y}^{\mfu,\alpha}(\theta^{\mfu,\alpha})\ge v_{*}(\theta^{\mfu,\alpha},X_{t,x}^{\alpha}(\theta^{\mfu,\alpha}))\;\;\as\;\forall\; \alpha \in \Ac^{t}.
$$
\end{theorem}

\proof The ingredients of the proof are essentially known, so we confine ourselves to a sketch. As $y>v(t,x)$, the definition of $v$ shows that there exists $\mfu \in \mfU^{t}$ satisfying
\be\label{eq: target reached proof gdp1}
Y_{t,x,y}^{\mfu,\alpha}(T)\ge g(X_{t,x}^{\alpha}(T))\;\;\as\;\forall\; \alpha \in \Ac^{t}.
\ee

\emph{Step 1.} We first consider the case where $\theta^{\mfu,\alpha}\equiv s \in [t,T]$ is a deterministic time independent of $\mfu,\alpha$. To be able to write processes as functionals of the canonical process, we pass to the raw filtration. More precisely, $\bar \F^{t}$ denotes the raw filtration generated by $(W_{s}-W_{t})_{t\le s\le T}$, extended trivially to $[0,T]$.
By \cite[Appendix~I, Lemma~7]{DellacherieMeyer.82}, we can find for each $\alpha \in \Ac^{t}$ an $\bar \F^{t}$-predictable process $\bar \mfu[\alpha]$ which is indistinguishable from  $\mfu[\alpha]$. The map $\alpha\in \Ac^{t}\mapsto \bar \mfu[\alpha]$ still satisfies the non-anticipativity condition~\reff{eq: cond non anticipativity} and therefore defines an element of $\mfU^{t}$. Moreover, \reff{eq: target reached proof gdp1} still holds if we replace $\mfu$ by $\bar \mfu$:
\be\label{eq: target reached proof gdp1 bis}
Y_{t,x,y}^{\bar \mfu,\alpha}(T)\ge g(X_{t,x}^{\alpha}(T))\;\;\as\;\forall\; \alpha \in \Ac^{t}.
\ee
We claim that it suffices to show that
\be\label{eq: claim for deterministic stopp times}
Y_{t,x,y}^{\mfu,\bar \alpha}(s) \ge      v_{*} (s,X_{t,x}^{\bar \alpha}(s))\;\;\as\;\mbox{ for all }\bar \alpha \in \bar \Ac^{t},
\ee
where $\bar \Ac^{t}$ is the set of all $\bar \F^{t}$-predictable processes with values in $A$.
Indeed, if $\alpha \in \Ac^{t}$, then by \cite[Appendix~I, Lemma~7]{DellacherieMeyer.82} we can find $\bar \alpha \in \bar \Ac^{t}$ such that $\alpha$ and $\bar \alpha$ are indistinguishable. In view of the non-anticipativity condition~\reff{eq: cond non anticipativity}, $\mfu[\alpha]$ and $\mfu[\bar \alpha]$ are also indistinguishable, and then \eqref{eq: claim for deterministic stopp times} implies that the same inequality holds for $\alpha$; that is, \eqref{eq: claim for deterministic stopp times} extends from $\bar \Ac^{t}$ to $\Ac^{t}$.

To prove~\eqref{eq: claim for deterministic stopp times}, fix $\bar \alpha\in \bar \Ac^{t}$.
For given $\omega\in \Omega$, we define
 $$
 \bar{\mfu}_{\omega}:(\tilde \omega, \tilde \alpha) \in \Omega \x \Ac^{s}\mapsto \bar \mfu[\bar \alpha(\omega)\oplus_{s}   \tilde \alpha](\omega \oplus_{s} (\tilde \omega -\tilde\omega_s +\omega_{s})),
$$
where we use the notation
$$
\gamma \oplus_{s} \gamma':=\gamma\1_{[0,s]} + \1_{(s,T]} \gamma'.
$$
We observe that $\bar{\mfu}_{\omega}\in \mfU^{s}$.
Using~\reff{eq: target reached proof gdp1 bis} and the flow property of the SDE~\eqref{eq: dynamics},
we can find a  nullset $N$ (depending on $\bar \alpha$) such that for all $\omega \notin N$ and all $\tilde{\alpha}\in\Ac^{s}$,
$$
Y_{s ,x'(\omega),y'(\omega)}^{\bar{\mfu}_{\omega},\tilde \alpha}(T) \ge g\big(X_{s ,x'(\omega)}^{\tilde \alpha}(T)\big) \;\;\as,
$$
where $x'(\omega):=X_{t,x}^{\bar \alpha}(s)(\omega)$ and $y'(\omega):=Y_{t,x,y}^{\bar \mfu,\bar \alpha}(s)(\omega)$.
By the definition of $v$, this means that $y'(\omega)\ge v(s,x'(\omega))$ or
$$
Y_{t,x,y}^{\bar \mfu,\bar \alpha}(s)(\omega) \ge v (s,X_{t,x}^{\bar \alpha}(s)(\omega))\;\;\forall \; \omega \notin N.
$$
Since $Y_{t,x,y}^{\bar \mfu,\bar \alpha}(s)=Y_{t,x,y}^{ \mfu,\bar \alpha}(s)$ a.s.~and $v\ge v_{*}$, this shows that~\eqref{eq: claim for deterministic stopp times} holds.

 \emph{Step 2.} To deduce the case of a general family $\{\theta^{\mfu,\alpha},(\mfu,\alpha)\in \mfU^{t}\x \Ac^{t}\}$ from Step~1, we approximate each $\theta^{\mfu,\alpha}$ from the right by stopping times with finitely many values and use the lower-semicontinuity of $v_{*}$ and the right-continuity of the paths of the state process. See e.g.~\cite[Section 2.3, Step 4]{BouchardMoreauNutz.12} for a very similar argument.
\ep
\\

We shall prove the supersolution property under the two subsequent conditions. A more general framework could be considered here (see e.g.~\cite{BoElTo09} or \cite{BouchardMoreauNutz.12}), but we shall anyway need these conditions for the second part of the dynamic programming principle below.
Given $(t,x,y,z,a)\in \D\x \R \x \R^{d}\x A$, define the set
$$
N(t,x,y,z,a):=\{u\in U: \sigma_{Y}(t,x,y,u,a)=z\}.
$$
The first condition is that $u\mapsto \sigma_{Y}(t,x,y,u,a)$ is invertible, and more precisely:

\begin{assumption}\label{ass: def hat u + regu} There exists a measurable map $\hat u:\D\x \R\x \R^{d} \x A\to U$  such that
$ N=\{\hat u\}$. Moreover, the map  $\hat u(\cdot,a)$ is continuous for each $a\in A$.
\end{assumption}

The second assumption is for the boundary condition at time $T$.

\begin{assumption}\label{ass: sequence Qna for T boundary}  Fix $a\in A$ and $(x,y)\in \R^d\x \R$. If there exist a sequence $(t_{n},x_{n},y_{n})\in \Di\x \R$ such that $(t_{n},x_{n},y_{n})\to (T,x,y)$
and a sequence $\mfu_{n}\in \mfU^{t_{n}}$ such that $Y_{t_{n},x_{n},y_{n}}^{\mfu_{n},a}(T)\ge g(X_{t_{n},x_{n}}^{a}(T))$ $\as$ for all $n\ge 1$, then $y\ge g(x)$.
\end{assumption}

\begin{remark} It follows from~\reff{eq: cond mu sigma} and~\reff{eq: ass v bounded g bounded lipschitz} that $g(X_{t_{n},x_{n}}^{a}(T)) \to g(x)$  a.s.\ for $n\to \infty$. If $U$ is bounded, then similarly $Y_{t_{n},x_{n},y_{n}}^{\mfu_{n},a}(T)\to y$ and we infer that Assumption~\ref{ass: sequence Qna for T boundary} holds. In the applications we have in mind, the assumption is satisfied even when $U$ is unbounded; see the proof of Corollary~\ref{cor: BSB} below.
\end{remark}

To state the supersolution property, let us define for $(t,x,y,q,p,M)\in \D\x \R\x \R\x \R^{d}\x \M^{d}$ the operators
\begin{eqnarray*}
  L^{a}(t,x,y,q,p,M) &:= &\mu^{\hat u}_{Y}(t,x,y,\sigma_{X}(t,x,a)p,a)-q-\mu_{X}(t,x,a)^{\top} p \\
  &&- \frac12[\sigma_{X}\sigma_{X}^{\top}(t,x,a)M]
\end{eqnarray*}
and
$$
  L:=\min_{a\in A}  L^{a},
$$
where
$$
\mu^{\hat u}_{Y}(t,x,y,z,a):=\mu_{Y}(t,x,y,\hat u(t,x,y,z,a),a),\;z\in \R^{d}.
$$

\begin{theorem}\label{thm: super sol property} Let Assumptions~\ref{ass: def hat u + regu} and~\ref{ass: sequence Qna for T boundary}  hold. Then, the function $v_{*}$ is a bounded viscosity supersolution of
\be\label{eq: supersol property}
\begin{array}{rcl}
  L(\cdot,\vp,\partial_{t}\vp,D\vp,D^{2}\vp)=0 &\mbox{on}& \Di\\[.3em]
 \vp-g=0&\mbox{on}& \DT.
 \end{array}
\ee
\end{theorem}

\proof
Now that GDP1 has already been established, the argument is similar to \cite[Section~5.1]{BoElTo09} and \cite[Theorem~3.4, Step~1]{BouchardMoreauNutz.12}. We sketch the proof for completeness.

\emph{Step 1.}  We start with the boundary condition, which is in fact an immediate consequence of Assumption~\ref{ass: sequence Qna for T boundary}. Let $(t_{n},x_{n})\to (T,x)$ be such that $v_{*}(t_{n},x_{n})\to v_{*}(T,x)$, for some $x \in \R^{d}$. Set $y_{n}:=v(t_{n},x_{n})+n^{-1}$. Then, we can find $\mfu_{n}\in \mfU^{t_{n}}$ such that  $Y_{t_{n},x_{n},y_{n}}^{\mfu_{n},a}(T)\ge g(X_{t_{n},x_{n}}^{a}(T))$ $\as$ for all $a\in A$. Sending $n\to \infty$ and using Assumption~\ref{ass: sequence Qna for T boundary} yields $v_{*}(T,x)\ge g(x)$ as desired.

\emph{Step 2.} We now prove that $v_{*}$ is a viscosity supersolution of~\reff{eq: supersol property} on $\Di$.
 Let $\vp$ be a smooth function and let $(t_{o},x_{o})\in \Di$ be such that
\be\label{eq: proof super sol min = 0}
\min_{\Di}{(\rm strict)} (v_{*}-\vp)= (v_{*}-\vp)(t_{o},x_{o})=0;
\ee
 we have to show that
$
L(\cdot,\vp,\partial_{t}\vp,D\vp,D^{2}\vp)(t_{o},x_{o})\geq 0$.
Suppose to the contrary that  we can find $a_{o}\in A$ such that
 $$
  L^{u_{o},a_{o}}(\cdot,\vp,\partial_{t}\vp,D\vp,D^{2}\vp)(t_{o},x_{o})<0,
 $$
 where $u_{o}:=\hat u(\cdot,\vp,\sigma_{X}(\cdot,a_{o})D\vp,a_{o})(t_{o},x_{o})$ and
 $$
  L^{u,a}(t,x,y,q,p,M) :=\mu_{Y}(t,x,y,u,a)-q-\mu_{X}(t,x,a)^{\top} p - \frac12[\sigma_{X}\sigma_{X}^{\top}(t,x,a)M].
 $$
The continuity of $\hat{u}$, cf.\ Assumption~\ref{ass: def hat u + regu}, implies that for all $\eps>0$, we can find $\delta>0$ such that     $|u-u_{o}|\le \eps$ whenever $u\in U$ is such that
$$
  |\sigma_{Y}(t,x,y,u,a_{o})-\sigma_{X}(t,x,a_o)D\vp(t,x)|\le \delta
$$
for some $(t,x,y)\in \Dc\x \R$ satisfying $|(t,x,y)-(t_{o},x_{o}, \vp(t_{o},x_{o}))|\le \delta$.
Recalling the regularity assumptions~\reff{eq: cond mu sigma} imposed on the coefficients of our controlled dynamics, this implies that we can find $\delta>0$ and an open neighborhood $B\subset \Di$ of $(t_{o},x_{o})$ such that
\be\label{eq: sub sol strict local proof super sol v*}
\begin{array}{c}
L^{u,a_{o}}(\cdot,y,\partial_{t}\vp,D\vp,D^{2}\vp)(t,x) \le 0\;\;\;

\forall\; (t,x)\in B  \;\mbox{ and }\;(y,u)\in \R\x U \; \mbox{ s.t. }
\\[.3em]
 |y-\vp(t,x)|\le \delta \;  \mbox{ and }\; |\sigma_{Y}(t,x,y,u,a_{o})-\sigma_{X}(t,x,a_{o})D\vp(t,x)|\le \delta.
\end{array}
\ee
We now fix $0<\eps<\delta$. It follows from~\reff{eq: proof super sol min = 0} that  we can  find   $(  t,  x,y)\in B \x \R$ satisfying
\be\label{eq: def eta et distance faible}
\vp(t,  x)+\eps>y>v(t,  x).
\ee
 Then, Theorem~\ref{thm: GDP1}
implies that  there exists $\mfu \in \mfU^{t}$ such that
\be\label{eq: proof super sol Y ge v*}
Y(s\wedge \theta)\ge v_{*}(s\wedge \theta,X(s\wedge \theta))\;\;\as\;\forall\;s\ge t,
\ee
where we abbreviate $(X,Y):=(X^{a_{o}}_{t,x},Y^{\mfu,a_{o}}_{t,x,y})$ and $\theta$ is the stopping time
$$
\theta:=\theta^{1}\wedge \theta^{2},
$$
$\theta^{1}$ is the first exit time of $(\cdot,X)$ from $B $ and $\theta^{2}$ is the first time $|Y -\vp(\cdot,X)|$ reaches the level $\delta$. Using~\reff{eq: proof super sol min = 0} again,  we observe that $v_{*}-\vp\ge 0$ on $B $ and  $(v_{*}-\vp)\ge \zeta$ on $\partial B $, for some $\zeta>0$. Then,~\reff{eq: proof super sol Y ge v*} and the definition of $\theta$ imply that
\begin{equation}\label{eq: proof sur sol Y ge vp}
\Delta(s):=Y(s\wedge\theta)-  \vp(s\wedge\theta,X(s\wedge\theta))\ge (\zeta\wedge \delta)\1_{\[\theta,T\]}(s)\;\;\as\;\forall\;s\ge t.
\end{equation}
 We now use~\reff{eq: sub sol strict local proof super sol v*} to obtain a contradiction to~\reff{eq: proof sur sol Y ge vp}. To this end, set
\begin{eqnarray*}
\lambda&:=&\sigma_{Y}(\cdot,X,Y,\mfu[a_{o}],a_{o})-\sigma_{X}(\cdot,X,a_{o})D\vp(\cdot,X),\\
\beta&:=& \left(L^{\mfu[a_{o}]_{\cdot},a_{o}}(\cdot,Y,\partial_{t}\vp,D\vp,D^{2}\vp)(\cdot,X)\right)  |\lambda|^{-2}\lambda \1_{\{|\lambda|\ge \delta\}}.
\end{eqnarray*}
It follows from the definition of $\theta$ and our regularity and relative growth conditions~\reff{eq: cond mu sigma} and~\reff{eq: cond mu over sigma}\footnote{This last condition is missing in \cite{BoElTo09}.} that $\beta$ is uniformly bounded on $\[t,\theta\]$.
This ensures that the (positive) exponential local martingale $M$ defined by the stochastic differential equation
$$
M(\cdot)=1-\int_{t}^{\cdot\wedge \theta} M(s) \beta_{s}^{\top}dW_{s}
$$
 is a true martingale. Recalling~\reff{eq: sub sol strict local proof super sol v*}, an application of It\^o's formula shows that $M\Delta$ is a local super-martingale. By the definition of $\theta\le \theta^{2}$, this process is bounded by the martingale $M\delta$, and is therefore a super-martingale.  In particular, \reff{eq: proof sur sol Y ge vp}   implies that
$$
y-\vp(t,x)=\Delta(t)\ge\Esp{M(\theta)\Delta(\theta)}\ge \Esp{M(\theta)( \zeta\wedge \delta)}=\zeta\wedge \delta.
$$
The required contradiction is obtained by choosing $\eps:=(\zeta\wedge \delta)/2$ in~\reff{eq: def eta et distance faible}.
\ep

\subsection{Second part of the dynamic programming principle and subsolution property}

We now turn to the second part of the geometric dynamic programming, GDP2, and thus to the main contribution of this paper. As already mentioned in the Introduction, we cannot rely on an abstract measurable selection argument as in \cite{SonerTouzi.02b}. Instead, we construct an almost optimal Markovian strategy that will play the role of a measurable selector in the proof of Theorem~\ref{thm: GDP2} below. This strategy is defined in~\reff{eq: def bar u verif gdp2} in terms of a  smooth supersolution $w$ of~\reff{eq: supersol property} having specific properties. The existence of $w$ will be proved in Section~\ref{sec: smoothing} by considering a family of controlled forward-backward stochastic differential equations and using the regularization technique of \cite{Krylov00}.
The arguments in Section~\ref{sec: smoothing} require existence and stability properties for the forward-backward stochastic differential equations, which hold under the following condition.
\begin{assumption}\label{ass: regu muYhatu}  {The map
$({t},x,y,z)\in {\D}\x\R\x\R^{d}\mapsto \mu_{Y}^{\hat u}(t,x,y,z,a)$ is {continuous and uniformly}   Lipschitz {in $(x,y,z)$}, uniformly in  ${a\in A}$, and  
$(y,z)\in \R\x\R^{d}\mapsto \mu_{Y}^{\hat u}(t,x,y,z,a)$ has linear growth, uniformly in  $(t,x,a)\in \D\x A$.}
\end{assumption}

We also assume that the comparison principle holds for~\reff{eq: supersol property}; see e.g.\ \cite{CrIsLi92} for sufficient conditions.

\begin{assumption}\label{ass: comparison}  Let $\underline w$ (resp.\ $\overline w$) be a lower-semicontinuous (resp.\ upper-semi\-continuous) bounded viscosity supersolution (resp.\ subsolution) of~\reff{eq: supersol property}. Then,  $\underline w\ge \overline w$ on $\D$.
\end{assumption}

As in \cite{Krylov00}, the regularization procedure also requires a concavity property for the operator $L^{a}$.
\begin{assumption}\label{ass: concavity} For all $(t,x,a)\in \D\x A$,
$
(y,p)\in \R\x \R^{d}\mapsto   L^{a}(t,x,y,0,p,0)
$
is concave.
\end{assumption}

We denote by  $C^{\infty}_{b}(\D)$ the set of bounded functions on $\D$ which have bounded derivatives of all orders. The following lemma summarizes the result of Section~\ref{sec: smoothing} in the present context.

\begin{lemma}\label{lem : smoothing supersol semi lin conv}  Let Assumptions~\ref{ass: def hat u + regu} to~\ref{ass: concavity} hold. Let $B\subset \D$ be a compact set and  let $\phi$ be a   continuous function such that $\phi\ge v_{*}+\eta$ on $B$, for some $\eta>0$. Then, there exists $w\in C^{\infty}_{b}(\D)$ such that
\begin{enumerate}
\item  $w$ is a classical supersolution of~\reff{eq: supersol property},
\item  $w \le \phi$ on $B$.
\end{enumerate}
\end{lemma}

\proof By Theorem~\ref{thm: super sol property} and Assumption~\ref{ass: comparison}, $v_{*}$ dominates any bounded subsolution of~\reff{eq: supersol property}. Thus, the lemma is a special case of Theorem~\ref{thm : smoothing supersol semi lin conv gene} below, applied with $\underline{w}:=v_{*}$. \ep \\

 We can now state  the main result of this section, the second part of the geometric dynamic programming principle.

\begin{theorem} [{\bf GDP2}]\label{thm: GDP2} Let Assumptions~\ref{ass: def hat u + regu} to~\ref{ass: concavity} hold.  Fix $(t,x,y)$ in $\D\x \R$, let $B\subset \D$ be a compact set containing $(t,x)$ and let $\theta^{\alpha}$ be the first exit time of $(\cdot, X^{\alpha}_{t,x})$ from $B$, for $\alpha \in  \Ac^{t}$.  Let  $\phi$ be a   continuous function such that $\phi\ge v_{*}+\eta$ on $\partial B\setminus \DT$ for some $\eta>0$ and suppose that there exists $\mfu_{o} \in \mfU^{t}$ such that
$$
Y_{t,x,y}^{\mfu_{o},\alpha}(\theta^{ \alpha})\ge \phi(\theta^{\alpha},X_{t,x}^{\alpha}(\theta^{ \alpha}))\1_{\{\theta^{\alpha}<T\}}
+g(X_{t,x}^{\alpha}(T))\1_{\{\theta^{\alpha}=T\}} \;\;\as\;\forall\; \alpha \in \Ac^{t}.
$$
Then, there exists $\mfu \in \mfU^{t}$ such that
$$
Y_{t,x,y}^{\mfu,\alpha}(T)\ge g(X_{t,x}^{\alpha}(T))\;\;\as\;\forall\; \alpha \in \Ac^{t},
$$
and in particular $y\ge v(t,x)$.
\end{theorem}

\proof Note that the lower-semicontinuity of $v_{*}$ ensures that $\phi\ge v_{*}+\eta$ on the closure of the bounded set $ \partial B\setminus \DT$. It then follows from Lemma~\ref{lem : smoothing supersol semi lin conv} (applied to the closure of $\partial B\setminus \DT$) that we can find a function  $w\in C^{\infty}_{b}$ which is a supersolution of
\reff{eq: supersol property} and satisfies  $w \le \phi$ on   $ \partial B\setminus \DT$. Next, we introduce $\mfu\in\mfU^{t}$ satisfying
\be\label{eq: def bar u verif gdp2}
\mfu[\alpha]=\mfu_{o}[\alpha]\1_{\[t,\theta^{\alpha} \[}+\1_{\[\theta^{\alpha},T\]}\hat u(\cdot, X_{t,x}^{\alpha},Y^{\mfu,\alpha}_{t,x,y},(\sigma_{X}(\cdot,\alpha)Dw)(\cdot, X_{t,x}^{\alpha}),\alpha)
\ee
for $\alpha\in \Ac^{t}$. To this end, let $Y$ be the unique strong solution of the equation
\begin{eqnarray*}
Y&=& Y^{\mfu_o,\alpha}_{t,x,y}(\theta^\alpha)+\int_{\theta^\alpha}^{\theta^\alpha\vee \cdot}\mu_{Y}^{\hat{u}}\big(\cdot,Y(s), \sigma_{X}(\cdot,\alpha_{s})Dw,\alpha_{s}\big)(s,X_{t,x}^{\alpha}(s)) ds
\\
&&+
\int_{\theta^\alpha}^{\theta^\alpha\vee \cdot} \sigma_{X}(s,X_{t,x}^{\alpha}(s),\alpha_{s})Dw(s,X_{t,x}^{\alpha}(s))dW_{s}
\end{eqnarray*}
which {is well-posed}  by Assumption~\ref{ass: regu muYhatu}. Then, if we define $\mfu$ by
\begin{eqnarray*}
\mfu[\alpha]:=\mfu_{o}[\alpha]\1_{\[t,\theta^{\alpha} \[}+\1_{\[\theta^{\alpha},T\]}\hat u(\cdot, X_{t,x}^{\alpha},Y,(\sigma_{X}(\cdot,\alpha)Dw)(\cdot, X_{t,x}^{\alpha}),\alpha),
\end{eqnarray*}
it follows via the definition of $\hat u$ in Assumption~\ref{ass: def hat u + regu} that
$$
  \sigma_{X}(s,X_{t,x}^{\alpha}(s),\alpha_{s})Dw(s,X_{t,x}^{\alpha}(s)) = \sigma_{Y}(s,X_{t,x}^{\alpha}(s),Y(s),\mfu[\alpha]_{s},\alpha_{s})
$$
for $s\geq \theta^{\alpha}$; that is,~\eqref{eq: def bar u verif gdp2} is satisfied.
%
%
%
%
 Note that $\mfu$ indeed belongs to the set $\mfU^{t}$ of non-anticipating strategies, due to the fact that   $\theta^{\alpha}=\theta^{\alpha'}$ on $\{\alpha{|_{[t,s]}}=\alpha'{|_{[t,s]}}\}\cap \{\theta^{\alpha}\le s\}$ for all $s\in [t,T]$ and $\alpha,\alpha'\in \Ac^{t}$.

Let us now fix $\alpha \in \Ac^{t}$.   Since $\phi\ge w$ on $\partial B\setminus \DT$, the definition of  the first exit time $\theta^{\alpha}$ implies that
$$
Y_{t,x,y}^{\mfu,\alpha}(\theta^{\alpha})\ge  w(\theta^{\alpha},X_{t,x}^{\alpha}(\theta^{\alpha}))\1_{\{\theta^{\alpha}<T\}}
+g(X_{t,x}^{\alpha}(T))\1_{\{\theta^{\alpha}=T\}}\;\;\as
$$
Applying It\^{o}'s formula to the smooth function $w$ and using that it is a supersolution of~\reff{eq: supersol property}   then leads to
$$
Y_{t,x,y}^{\mfu,\alpha}(T)\ge  w(T, X_{t,x}^{\alpha}(T))\1_{\{\theta^{\alpha}<T\}}+g(X_{t,x}^{\alpha}(T))\1_{\{\theta^{\alpha}=T\}}\ge g(X_{t,x}^{\alpha}(T))\;\;\as
$$
as claimed. Since $\alpha \in \Ac^{t}$ was arbitrary, this implies that $y\ge v(t,x)$.
\ep
\\

As a consequence of Theorem~\ref{thm: GDP2}, we can now prove that the upper-semi\-continuous envelope $v^{*}$ of $v$ is a viscosity subsolution of~\reff{eq: supersol property}.

\begin{theorem}\label{thm: sub sol property} Let Assumptions~\ref{ass: def hat u + regu} to~\ref{ass: concavity} hold. Then, the function $v^{*}$ is a bounded   viscosity subsolution of~\reff{eq: supersol property}.
\end{theorem}

\proof
We remark that due to the present (game) context, the proof in \cite{BoElTo09} cannot be reproduced per se; see Remark~\ref{rem: diff du to GDP2} for details. Instead, we first consider a case where $\mu_{Y}^{\hat u}$ is non-decreasing in $y$ and then treat the general case by reduction.

\emph{Step 1.}  Let $\mu_{Y}^{\hat u}$ be non-decreasing in its third variable $y$. We only prove the subsolution property on $\DT$; the subsolution property on $\Di$ follows from  similar arguments based on Theorem~\ref{thm: GDP2}; see \cite[Section 3]{BouchardMoreauNutz.12} and \cite[Section 5]{BoElTo09}. Let $x_{o}\in \R^{d}$ be such that
\be\label{eq: proof subsol max strict equal 0}
\max_{\D}{\rm (strict)}(v^{*}-\vp)=(v^{*}-\vp)(T,x_{o})=0
\ee
for some smooth function $\vp$, and suppose for contradiction that
\be\label{eq: proof subsol pde contradiction}
 \vp(T,x_{o})-g(x_{o}) =:2\kappa>0.
\ee
 Define $\tilde \vp(t,x):=\vp(t,x)+\sqrt{T-t}$ for $(t,x)\in \D$; then~\reff{eq: proof subsol pde contradiction} implies that there exists $\delta>0$ such that
\be\label{eq: contra 1}
& \tilde \vp-g >\kappa \;   \mbox{ on } B_{\delta}:= \{(t,x)\in \D : |(t,x)-(T,x_{o})|\le \delta\}.~~~&
\ee
   Moreover, the fact that   $\partial_{t} \tilde \vp(t,x)\to -\infty$ as $t\to T$ and the monotonicity assumption on $\mu_{Y}^{\hat u}$ imply  that, after possibly changing $\delta>0$,
\begin{equation}\label{eq: contra 2}
\inf_{a\in A} L^{a}(\cdot,y,\partial_{t} \tilde \vp ,D\tilde \vp ,D^{2}\tilde \vp )(t,x) \ge 0,
\;\;\forall\;
(t,x,y)\in B_{\delta}\x \R\;\mbox{ s.t.\ }\;   y\ge \tilde \vp(t,x)- \delta.
\end{equation}
Let
\be\label{eq: contra def zeta}
-\zeta:=\sup_{\partial B_{\delta} \setminus  \DT } (v^{*}-\tilde \vp)<0,
\ee
where the strict inequality follows from the fact that $(T,x_{o})$ achieves  a strict maximum of $v^{*}-\tilde \vp$.
   Fix $0<\eps< \delta\wedge \zeta\wedge \kappa$ and let $(t,x)\in B_{\delta}$ and $y\in \R$ be such that
  $$
  -\eps+\tilde \vp(t,x)<y <v(t,x);
  $$
see~\reff{eq: proof subsol max strict equal 0}.      Next, consider the strategy defined in a Markovian way by
  $$
  \alpha \in \Ac^{t} \mapsto \hat \mfu[\alpha]:=\hat u(\cdot, X_{t,x}^{\alpha},Y^{\hat \mfu,\alpha}_{t,x,y},(\sigma_{X}(\cdot, \alpha)D\tilde \vp)(\cdot, X_{t,x}^{\alpha}),\alpha).
  $$
  Without loss of generality, we can assume that $D\tilde \vp$ is bounded on $[t,T]\x \R^{d}$ and then the same arguments as in the proof of  Theorem~\ref{thm: GDP2} show that $\hat \mfu$ is well-defined as an element of $\mfU^{t}$.

  Given $\alpha\in \Ac^{t}$, let $\theta^{\alpha}$ be the first exit time of $(\cdot, X_{t,x}^{\alpha})$ from $B_{\delta}$.
  Then, by the definition of $\hat u$ and~\reff{eq: contra 2}, the process  $Y^{\hat \mfu,\alpha}_{t,x,y}-\tilde \vp(\cdot,X_{t,x}^{\alpha})$ is non-decreasing until it reaches $-\delta$ or time reaches $\theta^{\alpha}$. But since this process starts above $-\eps>-\delta$, it is in fact non-decreasing until  $\theta^{\alpha}$ and hence
  \begin{eqnarray*}
  Y^{\hat \mfu,\alpha}_{t,x,y}(\theta^{\alpha})&\ge & \tilde \vp(\theta^{\alpha},X_{t,x}^{\alpha}(\theta^{\alpha})) -\eps.
  \end{eqnarray*}
  Set $\eta:=\zeta-\eps>0$. Then, $\phi:=\tilde \vp-\eps\ge v_{*}+\eta$  on the closure of  $\partial B_{\delta}\setminus \DT$; see~\reff{eq: contra def zeta}. Moreover, it follows from~\reff{eq: contra 1} that $\tilde \vp-\eps \ge g+\kappa-\eps\ge g$  on  $B_{\delta}\cap \DT$; in particular,
  \begin{eqnarray*}
  Y^{\hat \mfu,\alpha}_{t,x,y}(\theta^{\alpha})&\ge & \phi(\theta^{\alpha},X_{t,x}^{\alpha}(\theta^{\alpha}))\1_{\{\theta^{\alpha}<T\}} +g(X_{t,x}^{\alpha}(T))\1_{\{\theta^{\alpha}=T\}}.
  \end{eqnarray*}
  Since $\alpha \in \Ac^{t}$ was arbitrary, Theorem~\ref{thm: GDP2} yields a contradiction to the fact that $y<v(t,x)$.

\emph{Step 2.} We now turn to the general case. Fix $\rho>0$ and define $\tilde Y_{t,x,y}^{\mfu,\alpha}$ as the strong solution of
$$
d\tilde Y(s)=\tilde \mu_{Y}(s,X_{t,x}^{\alpha}(s),\tilde Y(s),\mfu[\alpha]_{s}, \alpha_{s}) ds +\tilde \sigma_{Y}(s,X_{t,x}^{\alpha}(s),\tilde Y(s),\mfu[\alpha]_{s},\alpha_{s})dW_{s}
$$
with initial data $Y(t)=y$, where
\begin{eqnarray*}
\tilde \mu_{Y}(t,x,y,u,a)&:=& \rho y+e^{\rho t} \mu_{Y}(t,x,e^{-\rho t}y,u,a), \\
\tilde \sigma_{Y}(t,x,y,u,a)&:=&e^{\rho t} \sigma_{Y}(t,x,e^{-\rho t} y,u,a).
\end{eqnarray*}
Set $\tilde g:=e^{\rho T} g$ and define
$$
\tilde v(t,x):= \inf\{y\in \R: \exists \; \mfu \in \mfU^{t}\;\mbox{ s.t.}\; \tilde Y^{\mfu,\alpha}_{t,x,y}(T)\ge \tilde g(X^{\alpha}_{t,x}(T))\;\as\;\forall\;\alpha \in \Ac^{t}\}.
$$
Since   $\mu_{Y}^{\hat u}$ has linear growth  in its second argument $y$, see Assumption~\ref{ass: regu muYhatu}, one can choose $\rho>0$ so that
$$
\tilde \mu_{Y}^{\hat u}: (t,x,y,z,a) \mapsto  \rho y +  e^{\rho t}\mu_{Y}^{\hat u}(x,e^{-\rho t}y,e^{-\rho t}z,a)
$$
is non-decreasing in its $y$-variable. This means that these dynamics satisfy the monotonicity assumption used in Step~1 above; moreover, Assumptions~\ref{ass: def hat u + regu} to~\ref{ass: concavity} are also satisfied.  Hence, the upper-semicontinuous envelope  $\tilde v^{*}$ of $\tilde v$   is a viscosity subsolution of
$$
\begin{array}{rcl}
  \tilde L(\cdot,\vp,\partial_{t}\vp,D\vp,D^{2}\vp)=0 &\mbox{on}& \Di\\
 \vp-\tilde g=0&\mbox{on}& \DT,
 \end{array}
 $$
 where $\tilde L$ is defined like $L$ but with $\tilde \mu_{Y}$ and $\tilde \sigma_{Y}$ instead of $\mu_{Y}$ and $\sigma_{Y}$; that is,
 $$
 \tilde L:=\inf_{a\in A} \tilde L^{a}
 $$
and
 \begin{eqnarray*}
 \tilde L^{a} (t,x,y,q,p,M) &:=&  \rho y+ e^{\rho t}\mu^{\tilde u}_{Y}(t,x,e^{-\rho t}y,e^{-\rho t}\sigma_{X}(t,x,a)p,a) \\&&-q-\mu_{X}(t,x,a)^{\top} p - \frac12[\sigma_{X}\sigma_{X}^{\top}(t,x,a)M],
 \end{eqnarray*}
 where $\tilde u$ is defined like $\hat{u}$ but now in terms of $\tilde \sigma_{Y}$.
 Since $\tilde v(t,x)=e^{\rho t} v(t,x)$, this is equivalent to saying that $v^{*}$ is a viscosity subsolution of~\reff{eq: supersol property}.
\ep

\begin{remark}\label{rem: diff du to GDP2} In the proofs of \cite[Section 5]{BoElTo09},  the condition $|y-\tilde \vp(t,x)|\le  \delta$ was used instead of $y\ge \tilde \vp(t,x)- \delta$ as in~\reff{eq: contra 2}.
Correspondingly, $\theta^{\alpha}$ would then be the minimum of the first time when $(\cdot,X^{\alpha}_{t,x})$ reaches the boundary of $B_{\delta}$ and the first time when $|Y^{\hat \mfu,\alpha}_{t,x,y}-\tilde \vp(\cdot, X^{\alpha}_{t,x})|$ reaches $\delta$. In this case, similar arguments as above imply that
 \begin{eqnarray*}
  Y^{\hat \mfu,\alpha}_{t,x,y}(\theta^{\alpha})&\ge &  \phi(\theta^{\alpha},X_{t,x}^{\alpha}(\theta^{\alpha}))\1_{\{\theta^{\alpha}<T\}}  +g(X_{t,x}^{\alpha}(T))\1_{\{\theta^{\alpha}=T\}}
  \end{eqnarray*}
and in the context of \cite{BoElTo09}, this was enough to obtain a contradiction. However, it is not the case in our situation: this stopping time $\theta^{\alpha}$ is not an exit time of $(\cdot,X^{\alpha}_{t,x})$, which is a key condition in our Theorem~\ref{thm: GDP2}.
\end{remark}

\begin{remark}
It follows from Assumption~\ref{ass: comparison}, Theorem~\ref{thm: super sol property} and Theorem~\ref{thm: sub sol property}   that $v$ is continuous and that $v$ is the unique bounded (discontinuous) viscosity solution of~\reff{eq: supersol property}, whenever Assumptions~\ref{ass: def hat u + regu} to~\ref{ass: concavity} hold.
\end{remark}


\section{On the regularization of concave non-linear PDEs}\label{sec: smoothing}

In this section, we prove Lemma~\ref{lem : smoothing supersol semi lin conv} above. We consider a more general setting in order to isolate the result  from the particular context of the preceding section; the general version of Lemma~\ref{lem : smoothing supersol semi lin conv} is stated in Theorem~\ref{thm : smoothing supersol semi lin conv gene} below. Our result is very much in the spirit of \cite[Theorem 2.1]{Krylov00} (see also \cite{BarlesJakobsen.02}), which we extend to our setting. Consider the parabolic equation
 \be\label{eq: pde H}
H(\cdot,\vp,\partial_{t}\vp,D\vp,D^{2}\vp)=0\;\mbox{ on }\; \Di\;\mbox{ and }\; \vp=g \;\mbox{ on } \DT,
\ee
where, for  $(y,q,p,M)\in  \R\x \R\x \R^{d}\x \M^{d}$,
 \begin{align*}
 H(\cdot,y,q,p,M)
 :=\inf_{a\in A}\left(f(\cdot,y,\sigma_{X}(\cdot,a)p,a)-q-\mu_{X}(\cdot,a)^{\top}p -\frac12 {\rm Tr}[\sigma_{X}\sigma_{X}^{\top}(\cdot,a)M]\right),
 \end{align*}
 for some continuous function $f: \D\x \R \x  \R^{d}\x A\to\R$ such that
 \be\label{eq: ass f Lipschitz}
 \begin{array}{c}
(y,z) \mapsto f(t,x,y,z,a) \mbox{ is concave  {and has linear growth} } \\
{\mbox{ uniformly in  $(t,x,a)\in \D\x A$},}\\
 ({t},x,y,z) \mapsto f(t,x,y,z,a) \mbox{ is {continuous} and {uniformly} Lipschitz {in $(x,y,z)$},}\\
\mbox{ uniformly in ${a\in A}$.}
 \end{array}
  \ee
We continue to assume that the continuity and growth conditions~\reff{eq: cond mu sigma} hold for $\mu_X$ and $\sigma_X$, and that $g$ is bounded and Lipschitz continuous.

Our aim is to provide a smooth supersolution of~\reff{eq: pde H} which is controlled by a given viscosity supersolution of the same equation, in the sense of Theorem~\ref{thm : smoothing supersol semi lin conv gene} below. To this end, we first introduce a family of supersolutions of perturbations of~\reff{eq: pde H}. Namely, they correspond to the \emph{shaken coefficients} in the terminology of \cite{Krylov00}; i.e.\ to the operators $H_{\eps}$ defined for $(t,x,y,q,p,M)   \in \D\x \R\x \R\x \R^{d}\x \M^{d}$ by
  $$
 H_{\eps}(t,x,y,q,p,M)  :=\inf_{b\in B_{\eps}(0)}H((t,x)+b,y,q,p,M)\;,\;\;\eps>0,\;
 $$
 where $B_{\eps}(0)\subset \R^{d+1}$ denotes the closed ball of radius $\eps$ around the origin. This is also the first step in the analysis of \cite[Theorem 2.1]{Krylov00}, where $f$ does not depend on $y$ and $z$ and the result is obtained by considering a family of standard optimal control problems in the so-called Bolza form. We shall obtain the extension to our framework by considering instead controlled forward-backward stochastic differential equations.

 As in the previous section, we impose that the comparison principle holds for~\reff{eq: pde H}; this will ensure that the family of supersolutions associated to $H_{\eps}$ is controlled from above by the supersolutions of~\reff{eq: pde H}:

\begin{assumption}\label{ass: comp H general}   Let $\underline w$ (resp.\ $\overline w$) be a lower-semicontinuous (resp.\ upper-semi\-continuous) bounded viscosity supersolution (resp.\ subsolution) of~\reff{eq: pde H}. Then,  $\underline w\ge \overline w$ on $\D$.
\end{assumption}

\begin{proposition}\label{prop: existence of a continuous solution to the shaked pde}
  Let Assumption~\ref{ass: comp H general} hold. For all $\eps\in [0,1]$, there exists a  bounded and continuous  viscosity solution $w_{\eps}$  of
  \be\label{eq : pde shaked}
    H_{\eps}(\cdot,\vp, \partial_{t}\vp, D\vp,D^{2}\vp) =0 \mbox{ on  $\Di$}
  \ee
  with the following properties:
  \begin{enumerate}
  \item There exists  $c^{\eps}\in [0,T)$  such that
  \be\label{eq: varpi iota ge g + iota arround T}
  w_{\eps} \ge g+  \eps \mbox{ on $[T-c^{\eps},T]\x \R^{d}$ and $c^{\eps}>0$ if $\eps>0$.}
  \ee

  \item For any compact set $B\subset \D$, there exists $c_{\eps}^{B}$ satisfying  $c_{\eps}^{B}\to 0$ as $\eps\to 0$ and
  \be\label{eq: regu varpi in iota}
  w_{\eps} -c_{\eps}^{B}   \le w_{0}  \mbox{ on $B$.}
  \ee
  \end{enumerate}
\end{proposition}

\proof We shall construct the functions $w_{\eps}$ as  value functions associated to controlled forward-backward stochastic differential equations. Some of the subsequent technical arguments are known and will only be sketched; we focus on the points specific to our problem.

\emph{Step 1.}   In this step, we introduce a family of controlled forward-backward equations indexed by $\eps$. To this end, let $\Tc$ denote the set of stopping times and $\Db_{2}$ the collection of pairs $(\tau,\xi)$ such that $\tau\in \Tc$ and  $\xi \in L^{2}(\Fc_{\tau};\R^{d})$. Moreover, we denote by $\Sb_{2}$ the set of predictable processes $\psi$ such that $\sup_{t\le T} |\psi_{t}| \in L^{2}$ and by $\Hb_{2}$ the set of $\P\otimes dt$-square-integrable predictable processes.
Finally, we let $\Bc_{\eps}$ (resp.\ $\Ac$) denote the set of predictable processes with values in the ball $B_{\eps}(0)\subset \R^{d+1}$ (resp.~$A$). For each $\eps\in[0,1]$, given a control $\gamma:=(\alpha,\beta)\in \Ac\x \Bc_{\eps} =:\Gamma_{\eps}$ and an initial condition $\zeta:=(\tau,\xi)\in \Db_{2}$, we now consider the decoupled forward-backward stochastic differential equation
\begin{eqnarray*}
\bar X^{\zeta,\gamma}&=& \xi +\int_{\tau}^{\cdot} \mu_{X}((s,\bar X^{\zeta,\gamma}_{s})+\beta_{s},\alpha_{s})ds+\int_{\tau}^{\cdot} \sigma_{X}((s,\bar X^{\zeta,\gamma}_{s})+\beta_{s},\alpha_{s})dW_{s},\\
\bar Y^{\eps,\zeta,\gamma}&=&g_{\eps}(\bar X^{\zeta,\gamma}_{T}) -\int_{\cdot\vee \tau}^{T}f^{\zeta}_{\gamma}(s, \bar Y^{\eps,\zeta,\gamma}_{s},\bar Z^{\eps,\zeta,\gamma}_{s}   ) ds -
\int_{\cdot\vee \tau}^{T} \bar Z^{\eps,\zeta,\gamma}_{s}   dW_{s},
\end{eqnarray*}
where  $g_{\eps}:=g+2\eps$ and
$$
f^{\zeta}_{\gamma}(s,y,z):=f((s,\bar X^{\zeta,\gamma}_{s})+\beta_{s}, y ,z,\alpha_{s}  ) \mbox{ for } (s,y,z)\in [0,T]\x \R\x \R^{d}.
$$
Note that the forward equation does not depend on $\eps$. Moreover,~\reff{eq: ass f Lipschitz} implies that
\be\label{eq: f lipschitz}
\begin{array}{c}
 (y,z)\in \R\x \R^{d}\mapsto f^{\zeta}_{\gamma}(s,y,z)\mbox{ is  Lipschitz continuous  with linear growth} \\
 \mbox{ uniformly in $(s,\zeta,\gamma)\in [0,T]\x \Db_{2}\x\Gamma_{1}$. }
\end{array}
\ee
 Existence and uniqueness of the solution $(\bar X^{\zeta,\gamma},\bar Y^{\eps,\zeta,\gamma},\bar Z^{\eps,\zeta,\gamma} )\in \Sb_{2}\x \Sb_{2}\x \Hb_{2} $ are standard; see e.g.\ \cite{pardoux1998backward}.
 For $\gamma \in \Gamma_{\eps}$, $\zeta=(\tau,\xi)\in \Db_2$ and $\theta\in \Tc$ such that $\theta \ge \tau$ $\as$, we now define
\be\label{eq: def Yc iota}
\Yc^{\eps,\zeta,\gamma}_{\theta}:=\esssup_{\gamma'\in \Gamma_{\eps}} \bar Y^{\zeta,\gamma\oplus_{\theta}\gamma'}_{\theta},
\;\mbox{ where }\;
\gamma\oplus_{\theta} \gamma':=\gamma\1_{\[0,\theta\]}+\1_{\]\theta,T\]} \gamma'.
\ee

By a Girsanov change of measure argument as in \cite[Proposition~3.3]{BuLi08} (which is in turn an extension of \cite{Peng.97b}), it follows that
\be\label{eq: def wiota}
(t,x)\in \D\mapsto w_{\eps}(t,x):=\Yc^{\eps,(t,x),\gamma}_{t}
\ee
is a deterministic function on $\D$, independent of $\gamma \in \Gamma_{1}$.
In the remainder of the proof, we show that $w_{\eps}$ satisfies the requirements of the proposition.


\emph{Step 2.} We provide some estimates that will be used later on.
We first observe that~{\reff{eq: ass f Lipschitz}},  \cite[Theorem 1.5]{pardoux1998backward},~\reff{eq: cond mu sigma} and the Lipschitz continuity of $g$  imply that we can find $c>0$, independent of $\eps$, such that
\begin{align}
&\Esp{  |\bar Y^{\eps,\zeta_{1},\gamma\oplus_{\tau_{1}}\gamma'}_{\tau_{1}}-{\mathbb{E}_{\tau_{1}}}\left[\bar Y^{\eps,\zeta_{2},\gamma\oplus_{\tau_{2}}\gamma'}_{\tau_{2}}\right]|^{2} } {\le } c \Esp{|\xi_{1}-\xi_{2}|^{2}+|\tau_{1}-\tau_{2}| }\label{eq: a priori controls on BSDE}
\end{align}
for all $\zeta_{i}=(\tau_{i},\xi_{i})\in \Db_2$, $i=1,2${, such that $\tau_{1}\le \tau_{2}$}, and $\gamma,\gamma' \in \Gamma_{\eps}$.   As the family $\{\bar Y^{\eps,\zeta_{1},\gamma\oplus_{\tau_{1}}\gamma'}_{\tau_{1}}\}_{\gamma'\in \Gamma_{\eps}}$ is directed upward, it follows from  \cite[Proposition VI-1-1]{Neveu.75} that we can find a sequence $(\gamma'_{n})_{n}\subset \Gamma_{\eps}$ such that $\bar Y^{\eps,\zeta_{1},\gamma\oplus_{\tau_{1}}\gamma'_{n}}_{\tau_{1}}\uparrow \Yc^{\eps,\zeta_{1},\gamma}_{\tau_{1}}$ $\as$; see e.g.\ \cite[Lemma 2.4]{BouchardMoreauNutz.12} for a similar argument. Then
$$
\Yc^{\eps,\zeta_{1},\gamma}_{\tau_{1}}-{\mathbb{E}_{\tau_{1}}}\left[\Yc^{\eps,\zeta_{2},\gamma}_{\tau_{2}}\right]\le \liminf_{n} \left(\bar Y^{\eps,\zeta_{1},\gamma\oplus_{\tau_{1}}\gamma'_{n}}_{\tau_{1}}-{\mathbb{E}_{\tau_{1}}}\left[\bar Y^{\eps,\zeta_{2},\gamma\oplus_{\tau_{2}}\gamma'_{n}}_{\tau_{2}}\right] \right).
$$
Similarly, we can find   $(\gamma''_{n})_{n}\subset \Gamma_{\eps}$ such that
$$
{\mathbb{E}_{\tau_{1}}}\left[\Yc^{\eps,\zeta_{2},\gamma}_{\tau_{2}}\right]-\Yc^{\eps,\zeta_{1},\gamma}_{\tau_{1}}\le \liminf_{n} \left({\mathbb{E}_{\tau_{1}}}\left[\bar Y^{\eps,\zeta_{2},\gamma\oplus_{\tau_{2}}\gamma''_{n}}_{\tau_{2}}\right] -\bar Y^{\eps,\zeta_{1},\gamma\oplus_{\tau_{1}}\gamma''_{n}}_{\tau_{1}}\right).
$$
Moreover, as $g$ is bounded and $f$ satisfies~\reff{eq: f lipschitz}, we deduce from \cite[Theorem 2.2]{peng1993backward} that
$$
 \esssup_{\gamma'\in \Gamma_{\eps}}(| \bar Y^{\zeta_{2},\gamma'}_{\tau_{2}}|+|\bar Y^{\zeta_{1},\gamma'}_{\tau_{1}} | )\in L^{2}.
$$
Then, combining the above with the dominated convergence theorem and the a-priori estimate~\reff{eq: a priori controls on BSDE} yields that
\be\label{eq: L2 conti of Yc}
\Esp{  |\Yc^{\eps,\zeta_{1},\gamma}_{\tau_{1}}-{\mathbb{E}_{\tau_{1}}}\left[\Yc^{\eps,\zeta_{2},\gamma}_{\tau_{2}}\right]|^{2}}
\le c \Esp{|\xi_{1}-\xi_{2}|^{2}+|\tau_{1}-\tau_{2}|}.
\ee
{Applying this inequality   to $\zeta_{1}=(t,x) \in \D$ and $\zeta_{2}=(T,x)$, we see that  
\be \label{eq: inf bound for Y}
 | \Yc^{\eps,(t,x),\gamma}_{t}- g_{\eps}(x) |\le c\sqrt{T-t} \;\;\;\forall\; (t,x,\gamma)\in \D\x \Gamma_{\eps}.
 \ee}

\emph{Step 3.} The fact that  $w_{\eps}(z)=\Yc^{\eps,z,\gamma}_{t}$ for all $z=(t,x)\in \D$ ensures that $w_{\eps}(\zeta)=\Yc^{\eps,\zeta,\gamma}_{\tau}$ for simple random variables $\zeta=(\tau,\xi)\in \Db_{2}$. {The estimate~\reff{eq: L2 conti of Yc} shows that    the  function $w_{\eps}$ is continuous and a simple approximation argument implies that $w_{\eps}(\zeta)=\Yc^{\eps,\zeta,\gamma}_{\tau} $   for all   $\zeta=(\tau,\xi)\in \Db_{2} $ and $\gamma \in \Gamma_{\eps}$.} In particular, for $z=(t,x)\in \D$,
\be\label{eq: varpi = Yc}
w_{\eps}(\tau,\bar X^{z,\gamma}_{\tau})=\Yc^{\eps,z,\gamma}_{\tau} \mbox{ for all }  \tau \in \Tc_{t},
\ee
and $\Yc^{\eps,z,\gamma}$ admits a modification with continuous paths.

\emph{Step 4.}  We now sketch the proof of the dynamic programming principle for $\Yc^{\eps}$, and thus for $w_{\eps}$. Fix $z=(t,x)\in \D$, $\gamma \in \Gamma_{\eps}$ and let $t\le \tau\le \theta\le T$ be two stopping times.
Let $\gamma_{n}\in \Gamma_{\eps}$ be such that $\Yc^{\eps,z,\gamma}_{\theta}\le \bar Y^{z,\gamma\oplus_{\theta}\gamma_{n}}_{\theta}+n^{-1}$.
Then,  by the stability result    \cite[Theorem 1.5]{pardoux1998backward}, one can find $c>0$ and a sequence $(\delta_{n})_{n\ge 1}$ of random variables converging to $0$ a.s.\ such that
\begin{eqnarray*}
\Yc^{\eps,z,\gamma}_{\tau}\ge \bar Y^{\eps,z,\gamma\oplus_{\theta}\gamma_{n}}_{\tau}=\Ec^{f_{\gamma}^{z}}_{\tau,\theta}[\bar Y^{\eps,z,\gamma\oplus_{\theta}\gamma_{n}}_{\theta}]
\ge   \Ec^{f_{\gamma}^{ z}}_{\tau,\theta}[\Yc^{\eps,z,\gamma}_{\theta}]-c\delta_{n},
\end{eqnarray*}
where $\Ec^{f_{\gamma}^{z}}_{\tau,\theta}[\chi]$ denotes the value at time $\tau$ of the first component of the solution to the backward stochastic differential equation with driver $f_{\gamma}^{z}$ and terminal condition $\chi\in L^{2}(\Fc_{\theta})$ at time $\theta$.
 By sending $n\to \infty$, this implies that
 \begin{eqnarray*}
\Yc^{\eps,z,\gamma}_{\tau} \ge     \Ec^{f_{\gamma}^{z} }_{\tau,\theta}[\Yc^{\eps,z,\gamma}_{\theta}] .
\end{eqnarray*}
Conversely, it follows from    the comparison principle \cite[Theorem~1.6]{pardoux1998backward}  that
\begin{eqnarray*}
\Yc^{\eps,z,\gamma}_{\tau} = \esssup_{\gamma'\in \Gamma_{\eps}} \bar Y^{\eps,z,\gamma\oplus_{\tau}\gamma'}_{\tau}\le \esssup_{\gamma'\in \Gamma_{\eps}} \Ec^{f_{\gamma\oplus_{\tau}\gamma'}^{z}}_{\tau,\theta}[\Yc^{\eps,z,\gamma\oplus_{\tau}\gamma'}_{\theta}],
\end{eqnarray*}
which, combined with the above, implies
 \be\label{eq : dpp for Yc}
\Yc^{\eps,z,\gamma}_{\tau} =  \esssup_{\gamma'\in \Gamma_{\eps}} \Ec^{f_{\gamma\oplus_{\tau}\gamma'}^{z}}_{\tau,\theta}[\Yc^{\eps,z,\gamma\oplus_{\tau}\gamma'}_{\theta}].
\ee

\emph{Step 5.}   We already know from Step~3 that $w_{\eps}$ is continuous. Then, standard arguments based on   the identification~\reff{eq: varpi = Yc} and the dynamic programming principle~\reff{eq : dpp for Yc}   show that $w_{\eps}$ is a continuous viscosity solution of~\reff{eq : pde shaked}; see \cite{peng1992generalized} and~\cite{BuLi08} for details.

\emph{Step 6.} We conclude with the remaining estimates. The bound~\reff{eq: inf bound for Y} and the fact that $g$ is bounded imply that $(w_{\eps})_{\eps\in [0,1]}$ is uniformly bounded.
Moreover, the estimate~\reff{eq: inf bound for Y} implies~\reff{eq: varpi iota ge g + iota arround T};
 it remains to prove~\reff{eq: regu varpi in iota}. We first observe that the family $(w_{\eps})_{\eps\in [0,1]}$ is non-decreasing by the comparison principle for backward stochastic differential equations and therefore admits a limit $\bar w_{0}\ge w_{0}$ as $\eps \to 0$. The stability of viscosity solutions, see e.g.\ \cite{barles91},  combined with Assumption~\ref{ass: comp H general} ensures that $\bar w_{0}$ is a continuous  bounded viscosity solution of~\reff{eq: pde H}. Since $w_{0}$ is   a bounded continuous viscosity solution of the same  equation, Assumption~\ref{ass: comp H general} implies that $\lim_{\eps\downarrow 0} w_{\eps}=\bar w_{0}=w_{0}$. By Dini's theorem, the convergence is uniform on compact sets, which is~\reff{eq: regu varpi in iota}.
\ep

\begin{remark}
   Consider again the setting of Section~\ref{sec: stochastic target games} with Assumptions~\ref{ass: def hat u + regu} to~\ref{ass: concavity}.
   For $f:=\mu^{\hat u}_{Y}$, the equations~\eqref{eq: pde H} and~\eqref{eq: supersol property} coincide and hence $v=w_{0}$ by uniqueness, where $w_{0}$ is defined in~\reff{eq: def wiota}. Thus, the controlled forward-backward stochastic differential equation~\reff{eq: def Yc iota} with $\eps=0$ is a \emph{dual problem} for the stochastic target game.
\end{remark}

We now apply the smoothing technique used in \cite[Section 2]{Krylov00} to the functions $w_{\eps}$ obtained in Proposition~\ref{prop: existence of a continuous solution to the shaked pde} to construct a suitable smooth supersolution of~\reff{eq: pde H}.  The stability results~\reff{eq: varpi iota ge g + iota arround T} and~\reff{eq: regu varpi in iota} play an important role in ensuring that the correct boundary condition at time $T$ is satisfied and  that the upper bound in~(ii) below holds true.

 \begin{theorem}\label{thm : smoothing supersol semi lin conv gene} Let Assumption~\ref{ass: comp H general} hold and let $\wu$ be a function on $\D$ which dominates   any bounded viscosity subsolution of~\reff{eq: pde H}. Moreover, let $B\subset\D$ be a compact set and  $\phi$ be a   continuous function such that $\phi\ge \wu+\eta$ on $B$, for some $\eta>0$. Then,  there exists $w\in C^{\infty}_{b}(\D)$ such that
\begin{enumerate}
\item $w$ is a classical supersolution of~\reff{eq: pde H},
\item $w \le \phi$ on $B$.
\end{enumerate}
\end{theorem}

\proof  The proof is provided for the sake of completeness; we follow closely~\cite{Krylov00} and~\cite{ishii1995equivalence}. Throughout, any function $w$ on $\Dc$ is extended to $\R\x\R^d$ by $w(t,x)=w(0,x)$ for $t<0$ and $w(t,x)=w(T,x)$ for $t>T$.

\emph{Step 1.} We first construct a semi-concave function which is a.e.~a supersolution  of~\reff{eq: pde H}   in the interior of the parabolic domain. 
Let $w_{\eps}$ be as in Proposition~\ref{prop: existence of a continuous solution to the shaked pde}. For $k\ge 1$, consider the quadratic inf-convolution
$$
  w_{\eps}^{k}(z):=\inf_{z'\in \D} ( w_{\eps}(z' )+k|z-z'|^{2}),\;z\in \D.
$$
Since $ w_{\eps}$ is continuous and {bounded}, the infimum is achieved at some point $\hat z_{k}(z)$. Note that $ w_{\eps}\ge     w_{\eps}^{k}\ge - |w_{\eps}|_{\infty}$, where $|\cdot|_{\infty}$ denotes the sup-norm on $\D$. Hence,
$$
k|z-\hat z_{k}(z)|^{2}=        w_{\eps}^{k}(z)-w_{\eps}(\hat z_{k}(z))\le 2|w_{\eps}|_{\infty}=:l.
$$
It follows that
\be\label{eq: dist hat z k - z }
|z-\hat z_{k}(z)|^{2}\le l/k=:(\rho_{k})^{2}.
\ee
In particular, $\hat z_{k}(z)\to z$ as $k\to \infty$ and thus, using again $w_{\eps}\ge w_{\eps}^{k}$ and the continuity of $w_{\eps}$,
$$
 w_{\eps}(z)\le \liminf_{k\to \infty} ( w_{\eps}(\hat z_{k}(z))+k|z-\hat z_{k}(z)|^{2}) \le \limsup_{k\to \infty}     w_{\eps}^{k}(\hat z_{k}(z))\le  w_{\eps}(z).
$$
 Let $\vp$ be a smooth function on $\D$ and let $z\in \Di$ be such that
 $$
 \min_{\D}(   w_{\eps}^{k}-\vp)=(   w_{\eps}^{k}-\vp)(z)=0.
 $$
 Then, for any  $z'\in \D$,
\begin{eqnarray*}
 w_{\eps}(\hat z_{k}(z)+z'-z)+k|\hat z_{k}(z)-z|^{2}-\vp(z')
&\ge &   w_{\eps}^{k}(z')-\vp(z')\\
&\ge &    w_{\eps}^{k}(z)-\vp(z)\\
&=& w_{\eps}(\hat z_{k}(z))+k|\hat z_{k}(z)-z|^{2}-\vp(z).
\end{eqnarray*}
Hence, the minimum of $z'\in \D \mapsto  w_{\eps}(\hat z_{k}(z)+z'-z)-\vp(z')$ is achieved by $z'=z$ and therefore
$$
(\partial_{t}\vp,D\vp,D^{2}\vp)(z)\in \bar \Pc^{-} w_{\eps}(\hat z_{k}(z)),
$$
where  $\bar \Pc^{-} w_{\eps}(\hat z_{k}(z))$ denotes the closed parabolic subjet of $w_{\eps}$ at $\hat z_{k}(z)$; see e.g.~\cite{CrIsLi92}.
In view of Proposition~\ref{prop: existence of a continuous solution to the shaked pde}, this shows that $   w_{\eps}^{k}$ is a viscosity supersolution of
$$
H_{\eps}(\hat z_{k}(\cdot),\vp,\partial_{t}\vp,D\vp,D^{2}\vp)\ge 0 \mbox{ on } [\rho_{k},T-\rho_{k}]\x \R^{d}.
$$
Take $k$ large enough so that $\rho_{k}\le \eps/2$. Then,~\reff{eq: dist hat z k - z } yields that
$$
\{\hat z_{k}(z)+  b, \;b\in   B_{\eps}(0)\}\supset  \{z+  b, \;b\in   B_{\eps/2}(0)\}.
$$
This implies that
$   w_{\eps}^{k}$ is a viscosity supersolution of
$$
\min_{b\in B_{\eps/2}(0)} H(\cdot+b,\vp,\partial_{t}\vp,D\vp,D^{2}\vp)\ge 0 \mbox{ on } [\rho_{k},T-\rho_{k}]\x \R^{d}.
$$

\emph{Step  2.} We now argue as in \cite{ishii1995equivalence} to construct from the previous step a smooth supersolution in the interior of the parabolic domain. Since $w_{\eps}^{k}$ is semi-concave, there exist $D_{2}^{abs}w_{\eps}^{k} \in L^{1}(\D)$ and  a Lebesgue-singular negative Radon measure $ D_{2}^{sing} w_{\eps}^{k}$ on  $\D$ such that
$$
D^{2}w_{\eps}^{k}(dz)=D^{2}_{abs}w_{\eps}^{k}dz+D^{2}_{sing} w_{\eps}^{k}\mbox{ in the distribution sense}
$$
and
$$
(\partial_{t} w_{\eps}^{k},Dw_{\eps}^{k},D^{2}_{abs}w_{\eps}^{k})\in \bar \Pc^{-}w_{\eps}^{k} \mbox{ a.e. on } [\rho_{k},T-\rho_{k}]\x \R^{d};
$$
see \cite[Section 3]{jensen1988maximum}.
Hence, Step~1 implies that
\begin{equation}\label{eq: pde bar v k ae}
\min_{b\in B_{\eps/2}(0)} H(\cdot+b,w_{\eps}^{k}, \partial_{t} w_{\eps}^{k},D w_{\eps}^{k},D^{2}_{abs}w_{\eps}^{k})\ge 0 \mbox{ a.e. on } [\rho_{k},T-\rho_{k}]\x \R^{d}.
\end{equation}

 Next, we mollify $w_{\eps}^{k}$. Let $\psi\in C_b^{\infty}$ be a non-negative function with support $[-1,0]\x [-1,1]^{d}$ such that $\int_{\R^{d+1}} \psi(z)dz=1$,  and set $\psi_{\delta}(z)=\delta^{-d-1}\phi(z/\delta)$ for $z\in  \R^{d+1}$ and $\delta>0$. For any bounded measurable function $w$ on $\D$ and $(t,x)\in [\delta,T]\x \R^d$, let
\begin{eqnarray*}
w\star \psi_{\delta}(t,x)&:= &
\int_{\R\x\R^{d}} w(t', x')\psi_{\delta}(t'-t,x'-x )dt' dx'
\\
& =&
\int_{\R\x\R^{d}} w(t+t',x+x')\psi_{\delta}(t',x' )dt' dx'.
\end{eqnarray*}
Using~\reff{eq: pde bar v k ae} and the fact that $(y,q,p,M)\in \R\x\R\x \R^{d}\x \M^{d}\mapsto H(\cdot,y,q,p,M)$ is concave due to~\eqref{eq: ass f Lipschitz}, we obtain for   $\delta<\eps/2 $ that
\begin{eqnarray*}
H_{0}\big(\cdot, w_{\eps}^{k}\star \psi_{\delta}, \partial_{t} w_{\eps}^{k}\star \psi_{\delta} ,(D w_{\eps}^{k})\star \psi_{\delta} ,(D^{2}_{abs}w_{\eps}^{k})\star \psi_{\delta} \big)  (z_{o})\ge 0
\end{eqnarray*}
for all $z_{o}\in  [\rho_{k}+\delta,T-\rho_{k}]\x \R^{d}$.  Note that, since $D^{2}_{sing}w_{\eps}^{k}\le 0$ and $\psi_{\delta}\ge 0$,
\begin{eqnarray*}
\int_{\R\x\R^{d}} \psi_{\delta}D^{2}_{abs}w_{\eps}^{k}dz  &\ge& \int_{\R\x\R^{d}} \psi_{\delta}D^{2}_{abs}w_{\eps}^{k}dz+\int_{\R\x\R^{d}} \psi_{\delta}D^{2}_{sing}w_{\eps}^{k}(dz)\\
&=&  \int_{\R\x\R^{d}} w_{\eps}^{k} D^{2}\psi_{\delta}dz
\end{eqnarray*}
by integration-by-parts.
Since $L$ is parabolic, this  shows that, on $  [\rho_{k}+\delta,T-\rho_{k}]\x \R^{d}$,
\begin{eqnarray*}
 H_{0}\big(\cdot, w_{\eps}^{k}\star \psi_{\delta} , \partial_{t}(w_{\eps}^{k}\star \psi_{\delta} ),D ( w_{\eps}^{k} \star \psi_{\delta} ),D^{2}(w_{\eps}^{k}\star \psi_{\delta}) \big)  \ge 0.
\end{eqnarray*}
By sending $k\to \infty$,  we obtain  that  $w_{\eps,\delta}:=w_{\eps}\star \psi_{\delta}$ is a supersolution of~\reff{eq: pde H} on $(\delta,T)$. Moreover, $w_{\eps,\delta} \in C^{\infty}_{b}(\D)$ since $\psi\in C^{\infty}_{b}(\D)$ and $w_{\eps}$ is bounded by Proposition~\ref{prop: existence of a continuous solution to the shaked pde}.

\emph{Step  3.} We now choose $\eps$ and $\delta$ so that the supersolution property holds at~$T$ and~(ii) of the theorem is satisfied for $w=w_{\eps,\delta}$. We first use that $\phi\ge \wu+\eta$ and that $\wu$ dominates bounded subsolutions of~\reff{eq: pde H}; in particular, $\underline{w}\ge w_{0}$ by Proposition~\ref{prop: existence of a continuous solution to the shaked pde}. Choose $\eps$ such that the constant $c^B_\eps$ in~\reff{eq: regu varpi in iota} satisfies $c^B_\eps\leq \eta/2$, then
$$
  w_\eps \leq w_0 + c^B_\eps \leq \phi - \eta + c^B_\eps \leq \phi - \eta/2.
$$
As $w_\eps$ is continuous, we have $w_{\eps,\delta}\to w_{\eps}$ uniformly on the compact set $B$. Thus, for $\delta>0$ small enough with respect to $\eta$,
$$
w_{\eps,\delta}\le \phi \mbox{ on } B,
$$
which is assertion (ii). For the supersolution property at $T$, we first appeal to~\reff{eq: varpi iota ge g + iota arround T} and choose $\delta\in  (0, c^{\eps})$ such that we also have $w_{\eps,\delta}\ge g \star \psi_{\delta} + \eps$ on $\DT$. Since $g$ is uniformly continuous, we can choose $\delta>0$ small enough so that $g \star \psi_{\delta}\ge g-\eps$ on $\DT$, and therefore $w_{\eps,\delta}\ge g$ on the parabolic time boundary.

\emph{Step  4.} The fact that $w_{\eps,\delta}$ is only a supersolution of~\reff{eq: pde H} on $(\delta,T]$ and not on $[0,T]$ is not restrictive; we can follow the same arguments on $[-1,T]$ instead of $[0,T]$ and obtain that $w_{\eps,\delta}$ is   a supersolution of~\reff{eq: pde H} on $(-1+\delta,T] \supset [0,T]$ for $0<\delta<1$.
\ep

\begin{remark}\label{rem: open domain}
  In some applications, one might want to restrict the spatial domain to an open set $\Oc\subset \R^{d}$, for instance if the process $X^{\alpha}$ of Section~\ref{sec: stochastic target games} is defined as a stochastic exponential. Then, the equation~\reff{eq: pde H} should be naturally set on $\D:=[0,T]\x \Oc$. In this more general context, the arguments of Proposition~\ref{prop: existence of a continuous solution to the shaked pde} can be reproduced without modification whenever the comparison principle of Assumption~\ref{ass: comp H general} holds for $\D=[0,T]\x \Oc$. Moreover, the proof of Theorem~\ref{thm : smoothing supersol semi lin conv gene} can be generalized to this context under the following additional assumption: there exist $\delta_{\Oc}>0$ and $\psi_{\Oc}\in C^{\infty}_{b}(\R^{d}\x \R^{d};\R_+)$ such that, for all $x\in \R^d$,
  \begin{enumerate}
  \item $x+\delta x'\in \Oc$ for all   $\delta\in (0,\delta_{\Oc})$  and $ x' \in   \R^{d}$ such that $\psi_{\Oc}(x,x')\ne 0$,
  \item  $\int_{\R^{d}} \psi_{\Oc}(x,x')dx'=1$,
  \item the support of $\psi_{\Oc}(x,\cdot)$ is contained in $[-1,1]^{d}$.
  \end{enumerate}
  Indeed, it then suffices to replace the mollifier $\psi_{\delta}$ introduced in Step 2 of the proof of Theorem~\ref{thm : smoothing supersol semi lin conv gene} by
  $$
    \psi_{\delta}(t',x';x):=\delta^{-1-d}\phi(t'/\delta)\psi_{\Oc}(x,x'/\delta),\;\;(t',x',x)\in \R\x \R^{d}\x \R^{d},
  $$
  for some smooth function $\phi\ge0$ with  support $[-1,0]$ such that $\int_{\R} \phi(t)dt=1$, and to define the convolution operation by
  \begin{eqnarray*}
  w\star \psi_{\delta}(t,x)&:= &
  \int_{\R\x\R^{d}} w(t', x')\psi_{\delta}(t'-t,x'-x;x)dt' dx'
  \\
  & =&
  \int_{\R\x\R^{d}} w(t+t',x+x')\psi_{\delta}(t',x';x)dt' dx'.
  \end{eqnarray*}
  For  $\delta\le \delta_{\Oc}$, the condition (i) above ensures that $x'=x+\delta(x-x')/\delta \in \Oc$ whenever $x\in \Oc$ and $\psi_{\delta}(t'-t,x'-x;x)\ne 0$.

  Similarly, the proofs of Theorems~\ref{thm: GDP1} and~\ref{thm: super sol property} do not change if we replace  $[0,T]\x\R^{d}$ by  $[0,T]\x \Oc$. Under the additional assumption introduced above,   Lemma~\ref{lem : smoothing supersol semi lin conv} also extends to this case, and so do  the proofs of Theorems~\ref{thm: GDP2} and~\ref{thm: sub sol property}, whenever the comparison principle of Assumption~\ref{ass: comparison} holds for $\D=[0,T]\x \Oc$.
\end{remark}

\section{Super-hedging under model uncertainty}\label{sec: examples}

 In this section, we apply the above results to a super-hedging problem under model uncertainty, where the dynamics of the underlying stocks and the interest rates depend on a process $\alpha$ which is only known to take values in the given set $A$. We consider an investor who wants to hedge an option written on some underlying stocks whose $\log$-prices evolve according to
  $$
 X^{\alpha}_{t,x}=x+\int_{t}^{\cdot}  {\mu(s,X^{\alpha}_{t,x}(s),\alpha_{s})}ds +\int_{t}^{\cdot}  {\sigma(s,X^{\alpha}_{t,x}(s),\alpha_{s})}dW_{s},
 $$
 for some unknown process $\alpha \in \Ac^{t}$. Thus, there can be uncertainty both about the drift and the volatility. Moreover, the range of possible coefficients may be state-dependent. The investor also has a money market account at his disposal; however, the interest rates {$(r^{b}(s, X^{\alpha}_{t,x}(s),\alpha_{s}))_{s\le T}$ and $(r^{l}(s,X^{\alpha}_{t,x}(s),\alpha_{s}))_{s\le T}$} for borrowing and lending are different and also depend on the process $\alpha$.

A trading policy for the investor is a process $\nu\in \Uc^{t}$ taking values in $U=\R^{d}$. Recall that $\Uc^t$ consists of processes which are in $L^p(\P\x dt)$ for some $p\geq2$. In this section, we choose $p>2$. Each component $\nu^{i}$ corresponds to the monetary amount invested in the $i$-th stock. Given some initial capital $y\in \R$ at time $t$, the wealth process $Y^{\nu,\alpha}_{t,x,y}$ then evolves according to
\begin{eqnarray*}
Y^{\nu,\alpha}_{t,x,y}&=&y+\int_{t}^{\cdot} \nu_{s}^{\top} \{dX^{\alpha}_{t,x}{(s)}+\textstyle{\frac12}  {\gamma(s,X^{\alpha}_{t,x}(s),\alpha_{s})}ds\} \\
&&+\Int_{t}^{\cdot} \rho(s,{X^{\alpha}_{t,x}(s),}Y^{\nu,\alpha}_{t,x,y}(s),\nu_{s},\alpha_{s})ds\\
&=& y+{\int_{t}^{\cdot}\left( \nu_{s}^{\top} \{\mu+\textstyle{\frac12} \gamma\} +\rho\right) (s,X^{\alpha}_{t,x}(s),Y^{\nu,\alpha}_{t,x,y}(s),\nu_{s},\alpha_{s})ds}
\\
&&+ \int_{t}^{\cdot}  \nu_{s}^{\top}  {\sigma(s,X^{\alpha}_{t,x}(s), \alpha_{s})}dW_{s},
\end{eqnarray*}
where $\gamma$ is the vector containing the elements of the diagonal of the matrix $\sigma\sigma^\top$ and, with ${\bf 1}:=(1,\dots,1)^\top$,
$$
\rho(t,{x},y,u,a):=[y-u^{\top}{\bf 1}]^{+} r^{l}(t,{x},a)-[y-u^{\top}{\bf 1}]^{-} r^{b}(t,{x},a)
$$
for $(t,{x},y,u,a)\in [0,T]{\x\R^{d}}\x\R\x \R\x A$.

The process $\alpha$ may be interpreted as a directly observable risk factor such as the federal funds rate or an indicator for the state of the economy. In this case it clearly makes sense that the investor's actions may depend on $\alpha$ in  a non-anticipative way. Alternatively, as the stocks $X^{\alpha}_{t,x}$ are observable, so are the volatilities $\sigma(s,X^{\alpha}_{t,x}(s),\alpha_{s})$. Therefore, if $a\mapsto \sigma(s,x,a)$ is invertible, the process $\alpha$ is automatically observable. We thus allow the investor to choose a non-anticipative strategy $\mfu\in \mfU^{t}$. As a result, the
super-hedging price is defined by
$$
v(t,x):=\inf\{y\in \R:~\exists\; \mfu\in \mfU^{t}\;\mbox{ s.t.\ }\; Y^{\mfu,\alpha}_{t,x,y}(T)\ge g(X_{t,x}^{\alpha}(T))\;\as\;\forall\;\alpha \in \Ac^{t}\}
$$
or the smallest price which allows to super-hedge for all $\alpha \in \Ac^{t}$.

We henceforth assume that $g$ is bounded and Lipschitz continuous, that $\mu$ and~$\sigma$ are continuous functions, that
\be\label{eq: ass dyna example}
\begin{array}{c}
 (t,x,a)\in [0,T]\x\R^d\x A\mapsto (\mu,\sigma,{r^{l},r^{b}})(t,x,a)\;\mbox{ is {continuous}, bounded}, \\
 \mbox{{continuous in $(t,x)$} and Lipschitz continuous in $x$, {both} uniformly in ${a}$}, \\
 \mbox{and that }\sigma(t,x,a) \mbox{ is invertible for all }(t,x,a)\in [0,T]\x\R^d\x A.
 \end{array}
\ee
We also assume that 
$$
r^{b}\ge r^{l}
$$
and that
\be
    &\lambda^{b}:=\sigma^{-1}(\mu +\textstyle{\frac12}\gamma-r^b)\;\;\mbox{ and }\;\; \lambda^{l}:=\sigma^{-1}(\mu+\textstyle{\frac12}\gamma-r^l)\quad \mbox{are bounded,}&\label{eq: ass change measure rb rl}\label{eq: ass change measure rb rl}\\
    &{(t,x,a)\in [0,T]\x \R^{d}\x A\mapsto (\lambda^{b},\lambda^{l})(t,x,a) \;\mbox{ does not depend on  $x$.}}&\label{eq: ass indep x}
\ee
We observe that {the assumption  \reff{eq: ass change measure rb rl} is in fact a no-arbitrage condition}. In this setting, we have the following characterization of the super-hedging price.

\begin{corollary}\label{cor: BSB}
The function $v$ is continuous and the unique bounded viscosity solution of
\begin{eqnarray*}
\begin{array}{rcl}
\inf\limits_{a\in A}\left[\rho(\cdot,\vp,D\vp,a)+ \textstyle{\frac12}\gamma(\cdot,a)D\vp  -\Lc^{a}\vp\right]=0 &\mbox{ on } & [0,T)\x \R^{d}\\
\vp(T,\cdot)=g&\mbox{ on } &\R^{d},
\end{array}
\end{eqnarray*}
where
$$
\Lc^{a}\vp(t,x):=\partial_{t}\vp(t,x)
+\frac12{\rm Tr}[\sigma\sigma^{\top}(t,x,a) D^{2}\vp(t,x)].
$$
\end{corollary}

\proof The present model satisfies the conditions~\reff{eq: cond mu sigma} and~\reff{eq: cond mu over sigma} of Section~\ref{sec: stochastic target games}. Moreover, Assumption~\ref{ass: def hat u + regu} holds with
$$
\hat u(t,x,y,z,a)=(\sigma^{-1})^\top(t,x,a)z.
$$
To verify Assumption~\ref{ass: sequence Qna for T boundary}, consider $a\in A$ and $\mfu\in \mfU$. Then It\^o's formula, our choice $p>2$ and~\reff{eq: ass change measure rb rl} imply that the process $(e^{-\int_{t}^{\cdot} r^{b}(s,a)ds}$ $Y^{\mfu,a}_{t,x,y})$ is a super-martingale under the probability measure $\Q_{t,x}^{a}$  defined via the Dol\'eans-Dade exponential
$$
  \frac{d\Q_{t,x}^{a}}{d\P} = \Ec_T\left(-\int_{t}^{\cdot} \lambda^{b}(s, X_{t,x}^{a}(s), a) dW_{s} \right).
$$
In particular, if $\Di\x \R\ni (t_{n},x_{n},y_{n})\to (T,x,y)\in \DT\x \R$ and $(\mfu_{n})_{n\ge 1}$ are such that $\mfu_{n}\in \mfU^{t_{n}}$ and $ Y_{t_{n},x_{n},y_{n}}^{\mfu_{n},a}(T)\ge g(X_{t_{n},x_{n}}^{a}(T))$ $\as$, then $y_{n}\ge \E_{\Q_{t_{n},x_{n}}^{a}}[g(X_{t_{n},x_{n}}^{a}(T))]$ for all $n\ge 1$, while~\reff{eq: ass dyna example} and the continuity of the bounded function $g$ ensure that $\E_{\Q_{t_{n},x_{n}}^{a}}[g(X_{t_{n},x_{n}}^{a}(T))]\to g(x)$ as $n\to \infty$. Since $y_{n}\to y$, this shows that $y\ge g(x)$, so Assumption~\ref{ass: sequence Qna for T boundary} holds.

As $g$ is bounded, the same super-martingale property implies that $v$ is bounded from below. Moreover, as $r^l$ is bounded, $v$ is also bounded from above, so we have~\reff{eq: ass v bounded g bounded lipschitz}.
To verify Assumption~\ref{ass: regu muYhatu}, we use the definitions of $\hat{u}$ and $\rho$ to obtain that
\begin{eqnarray*}
 \mu^{\hat u}_{Y}(t,x,y,z,a):=\rho(t,{x,}y,\hat u(t,x,y,z,a),a)+\hat u(t,x,y,z,a)^{\top}\{\mu(t,x,a)+\textstyle{\frac12}\gamma(t,x,a)\}
\end{eqnarray*}
satisfies
\begin{eqnarray*}
| \mu^{\hat u}_{Y}(t,x,y,z,a)|&\le& |y|(|r^{b}|\vee |r^{l}|)(t,{x,}a)+|z|(|\lambda^{l}|\vee|\lambda^{b}|)(t,x,a)
\end{eqnarray*}
and
\begin{eqnarray*}
| \mu^{\hat u}_{Y}(t,x,y,z,a)- \mu^{\hat u}_{Y}(t,x,y',z',a)|&\le&|y-y'|(|r^{b}|\vee |r^{l}|)(t,{x,}a)\\
&&+|z-z'|(|\lambda^{b}|\vee |\lambda^{l}|)(t,x,a)
\end{eqnarray*}
for all $(t,x,a)\in [0,T]\x\R^d\x A$ and $(y,z),(y',z')\in \R\x \R^{d}$.  {The Lipschitz continuity in the $x$-variable follows immediately from \reff{eq: ass indep x}.} Assumption~\ref{ass: concavity} is satisfied\footnote{We may observe here that this concavity condition is satisfied in a very natural way in the context of mathematical finance, where wealth dynamics are typically either linear or exhibit concave nonlinearities. The concavity is related to the absence of arbitrage opportunities; in the present example, the fact that $r^b \ge r^l$.}  due to $r^{l}\le r^{b}$. Finally, Assumption~\ref{ass: comparison} follows from standard arguments under our Lipschitz continuity conditions; cf.\ \cite{CrIsLi92}. \ep

\bibliographystyle{plain}

\end{document}